\documentclass[12pt,a4paper]{amsart}

\usepackage[utf8]{inputenc}
\usepackage[margin=1in]{geometry}
\usepackage{hyperref}
\usepackage{amsmath}
\usepackage{amsfonts}
\usepackage{amsthm}
\usepackage{amssymb}
\usepackage{xcolor}
\usepackage{graphicx}
\usepackage{mathtools}
\usepackage{tikz-cd} 
\usepackage{bbm} 
\usepackage{enumerate} 
\usepackage{stmaryrd} 
\usepackage{subcaption} 
\usepackage{microtype} 

\hypersetup{linktocpage}

\newcounter{clm_cntr}
\newcounter{pf_clm_cntr}

\theoremstyle{plain}
\newtheorem{thrm}{Theorem}[section]
\def\bthm{\begin{thrm}}
\def\ethm{\end{thrm}}

\newtheorem{theoremx}{Theorem}

\newtheorem{prop}[thrm]{Proposition}
\def\bprop{\begin{prop}}
\def\eprop{\end{prop}}

\newtheorem{ques}[thrm]{Question}
\def\bques{\begin{ques}}
\def\eques{\end{ques}}

\newtheorem{cjec}[thrm]{Conjecture}
\def\bcjec{\begin{cjec}}
\def\ecjec{\end{cjec}}

\newtheorem{cor}[thrm]{Corollary}
\def\bcor{\begin{cor}}
\def\ecor{\end{cor}}

\newtheorem{fact}[thrm]{Fact}
\def\bfact{\begin{fact}}
\def\efact{\end{fact}}

\newtheorem{lem}[thrm]{Lemma}
\def\blem{\begin{lem}}
\def\elem{\end{lem}}

\theoremstyle{definition}
\newtheorem{defn}[thrm]{Definition}
\def\bdefn{\begin{defn}}
\def\edefn{\end{defn}}

\newtheorem*{conc}{Conclusion}
\def\bconc{\begin{conc}}
\def\econc{\end{conc}}

\newtheorem{alg}[thrm]{Algorithm}
\def\balg{\begin{alg}}
\def\ealg{\end{alg}}

\def\bproof{\begin{proof}}
\def\eproof{\end{proof}}

\theoremstyle{remark}

\newtheorem{rem}[thrm]{Remark}
\def\brem{\begin{rem}}
\def\erem{\end{rem}}

\newtheorem{ex}[thrm]{Example}
\def\bex{\begin{ex}}
\def\eex{\end{ex}}

\newtheorem{exs}[thrm]{Examples}
\def\bexs{\begin{exs}}
\def\eexs{\end{exs}}

\newtheorem{claim}[clm_cntr]{Claim}

\newtheorem*{pfsketch}{Sketch of Proof}
\newtheorem{pfclaim}[pf_clm_cntr]{Proof of Claim}
\newtheorem{obs}{Observation}
\def\bobs{\begin{obs}}
\def\eobs{\end{obs}}	

\let\ra=\rightarrow

\let\lra=\leftrightarrow

\let\bra=\mapsto

\let\tensor=\otimes

\let\setminus=\smallsetminus

\def\A{\mathcal{A}}
\def\H{\mathcal{H}}
\def\N{\mathbb{N}}
\def\Z{\mathbb{Z}}
\def\R{\mathbb{R}}
\def\C{\mathbb{C}}
\def\Q{\mathbb{Q}}

\def\K{\mathbb{K}}
\def\k{\mathbbm{k}}

\def\m{\mathfrak{m}}

\DeclareMathOperator{\hgt}{ht}

\DeclareMathOperator{\fchar}{char}

\DeclareMathOperator{\ft}{ft}
\DeclareMathOperator{\lct}{lct}
\DeclareMathOperator{\init}{in}
\DeclareMathOperator{\even}{even}
\DeclareMathOperator{\odd}{odd}

\DeclareMathOperator{\cut}{cut}
\DeclareMathOperator*{\maximize}{maximize}
\DeclareMathOperator{\subject}{subject \, to:}
\DeclareMathOperator*{\minimize}{minimize}
\DeclareMathOperator{\LP}{LP}
\DeclareMathOperator{\PP}{PP}
\DeclareMathOperator{\LPd}{LP^{\ast}}
\DeclareMathOperator{\PPd}{PP^{\ast}}
\DeclareMathOperator{\tr}{T}
\DeclareMathOperator{\bi}{bi}
\DeclareMathOperator{\MM}{MM}

\newcommand{\norm}[1]{|| #1 ||}

\newcommand{\card}[1]{\lvert #1 \rvert}

\newcommand{\ul}[1]{\underline{#1}}

\newcommand{\wh}[1]{\widehat{#1}}

\newcommand{\mbf}[1]{\mathbf{{#1}}}
\renewcommand{\subset}{\subseteq}

\AtBeginDocument{%
   \def\MR#1{}
}

\title{Invariants of Binomial Edge Ideals via Linear Programs}

\begin{document}
\author{Adam LaClair}
\thanks{The author was partially supported by NSF grant DMS-2100288 and by Simons Foundation Collaboration Grant for Mathematicians \#580839}
\address{Adam LaClair \\ Department of Mathematics \\ Purdue University \\
West Lafayette \\ IN 47907 \\ USA} 
\email{alaclair@purdue.edu}

\begin{abstract}
We associate to every graph a linear program for packings of vertex disjoint paths. We show that the optimal primal and dual values of the corresponding integer program are the binomial grade and height of the binomial edge ideal of the graph. We deduce from this a new combinatorial characterization of graphs of K\"onig type and use it to show that all trees are of K\"onig type.

The log canonical threshold and the F-threshold are important invariants associated to the singularities of a variety in characteristic $0$ and characteristic $p$. We show that the optimal value of the linear program (computed over the rationals) agrees with both the F-threshold and the log canonical threshold of the binomial edge ideal if the graph is a block graph or of K\"onig type. We conjecture that this linear program computes the log canonical threshold of the binomial edge ideal of any graph.

Our results resemble theorems on monomial ideals arising from hypergraphs due to Howald and others.
\end{abstract}

\maketitle
\setcounter{tocdepth}{1}
\tableofcontents

\section{Introduction}
The theory of matchings of a graph has a long and rich history within graph theory with important applications in industry and also to other branches of mathematics. There are classical theorems relating the combinatorics of hypergraphs, the linear program realizing a maximal matching in a hypergraph, and algebraic invariants of the edge ideal associated to a hypergraph. For a hypergraph $\H$ let $\MM_{\H,\Z}$ (respectively $\MM_{\H,\Q}$) denote the optimal value of the linear program associated to $\H$ realizing a maximal matching of $\H$ over $\Z$ (respectively $\Q$), and $\MM_{\H,\Z}^{\ast}$ denote the optimal value of the dual linear program over $\Z$. For background on hypergraphs and the edge ideal of a hypergraph see \cite{morey2012edge}. For background on (fractional) maximal matchings of a hypergraph or the associated linear program $\MM_{\H,-}$ see \cite{lovaz1986matching}. 
\begin{theoremx}[classical]
\label{thm_classical}
\label{thm_monomial_ideal_matching_combinatorics} 
Let $\H$ be a hypergraph, $I(\H)$ the associated edge ideal, and $\MM$ as above. 
    \begin{enumerate}
        \item The monomial grade of $I(\H)$ is equal to $\MM_{\H,\Z}$, and also to the size of a maximal matching of $\H$. (folklore) \label{item_thm_classical_1}
        \item The following quanitites are equal:
        \begin{itemize}
            \item For any prime $p$ the F-threshold of $I(\H)$,
            \item The log canonical threshold of $I(\H)$,
            \item $\MM_{\H,\Q}$,
            \item The size of a maximal fractional matching of $\H$.
        \end{itemize}
        This was proven by \cite{howald2001multiplier}, but see also \cite{hernandez2016f}, \cite{pagi2018enhanced}.\label{item_thm_classical_2}
        \item The height of $I(\H)$ is equal to $\MM^{\ast}_{\H,\Z}$, and also to the size of a minimal vertex cover of $\H$ (see for instance {\cite[Lemma 9.1.4]{herzog2011monomial}} and \cite{lovaz1986matching}). \label{item_thm_classical_3}
    \end{enumerate}
\end{theoremx}
Theorem \ref{thm_classical} can be formulated for an arbitrary monomial ideal if one is willing to sacrifice combinatorial notions (see \cite{howald2001multiplier}, \cite{hernandez2016f}, \cite{pagi2018enhanced}, \cite[Lemma 9.1.4]{herzog2011monomial}). 

Another classical theorem in matching theory is K\"onig's theorem.
\begin{theoremx}[K\"onig]
\label{thm_konig}
Let $G$ be a bipartite graph, then the size of a maximal matching is equal to the size of a minimal vertex cover.
\end{theoremx}
Finally, we recall results on the coordinates of the vertices of the fractional matching polytope.
\begin{theoremx}
\label{thm_vertices_matching_poly}
For a graph $G$ each coordinate of a vertex of the fractional matching polytope is either $0$, $1$, or $\frac{1}{2}$. (Balinski, see {\cite[Theorem 7.5.1]{lovaz1986matching}}. If $G$ is a bipartite graph, then each coordinate of the vertex of the fractional matching polytope is either $0$ or $1$. (see {\cite[Theorem 7.1.2]{lovaz1986matching}}.
\end{theoremx}
In this paper we introduce a new linear program associated to a graph. We prove analogues of the above statements relating the linear program, the combinatorics of the graph, and algebraic properties of the binomial edge ideal.

The class of binomial edge ideals was introduced in \cite{herzog2010binomial} and independently in \cite{ohtani2011graphs}. They are defined as 

\bdefn[\cite{herzog2010binomial}]
\label{defn_bin_edge_ideal}
Let $G = (V,E)$ be a finite simple graph with vertex set $V = \{1,\ldots,n\}$ and edge set $E$. Fix a field $\K$. Consider the polynomial ring $R := \K[X_{1},\ldots,X_{n},Y_{1},\ldots,Y_{n}]$, and for each edge $\{i,j\} \in E$ with $i < j$ define $f_{ij} := X_{i}Y_{j} - X_{j}Y_{i} \in R$. Define the binomial edge ideal of $G$, denoted $J_{G}$, to be the ideal
\begin{equation}
\label{eqn_binomial_edge_ideal}
J_{G} := ( \{ f_{ij} \mid \{i,j\} \in E \} ).
\end{equation}
\edefn
Numerous papers have studied algebraic invariants of a binomial edge ideal in terms of combinatorial properties of the underlying graph; see \cite{saeedi2016binomial} for a survey.

In \cite{takagi2004fpure}, the authors introduced the notion of F-pure threshold which is a characteristic $p$-invariant having properties analogous to log canonical thresholds. In our situation the F-pure threshold is equivalent to the F-threshold introduced in \cite{mustata2005FThresholds}. Throughout this paper $p$ will always denote a positive prime integer. When $S = \Z[X_{1},\ldots,X_{n}]$ and $I$ is an $S$-ideal, we denote by $I_{\C}$, respectively $I_{p}$, the image of $I$ in $S \tensor_{\Z} \C$, respectively $S \tensor_{\Z} \Z/p\Z$. We denote by $\lct(I)$ the log canonical threshold of $I_{\C}$ around the origin which is an important invariant appearing in birational geometry (see \cite{lazarsfeld2017positivity}) and by $\ft(I_{p})$ the F-threshold of $I_{p}$. Throughout this paper we will often write $\ft(J_{G})$ instead of $\ft(J_{G,p})$ since our computations will be independent of characteristic. It is known that:
 
\bthm[{\cite[Theorem 3.4]{mustata2005FThresholds}}]
\label{thm_ft_leq_lct_and_lct_equals_lim_ft}
In the above setup,
\begin{equation}
\ft(I_{p}) \leq \lct(I)
\end{equation}
for all $p \gg 0$ and
\begin{equation}
\lct(I) = \lim_{p \ra \infty} \ft(I_{p}).
\end{equation}
\ethm
In practice it is difficult to give explicit formulas for the log canonical threshold or F-threshold of an ideal, apart from special cases. Relevant for us are the following computations. The F-threshold and log canonical threshold of monomial ideals can be computed via a linear program (see \cite{howald2001multiplier}, \cite{hernandez2016f}, \cite{pagi2018enhanced}). Shibuta and Takagi \cite{shibuta2009logcanonical} used this linear program to compute the log canonical threshold of certain classes of binomial ideals as the optimal value of a linear programming problem over $\Q$. Hern\'{a}ndez \cite{hernandez2014Fpure} expands on their ideas to show how the optimal value of this linear program over $\Q$ can be used to compute the F-threshold of a binomial hypersurface. In \cite{miller2014fptDeterminantal}, the authors compute the F-threshold of determinantal ideals. This computation was later generalized by Henriques and Varbaro \cite{henriques2016test} where they compute the multiplier and test ideals of GL-equivariant ideals, and as a consequence they obtain a formula for the F-threshold of GL-equivariant ideals. In \cite{blanco2018procedure}, the authors give an algorithm to compute the log canonical threshold of any binomial ideal. Velez \cite{velez2019f} computes the F-threshold of the ideal of adjacent maximal minors which in the case of a $2\times n$ matrix corresponds to the binomial edge ideal of a path. Gonz{\'a}lez-Mart{\'\i}nez computes the F-threshold of the graded maximal ideal in the quotient ring $R/J_{G}$ for any graph $G$ and uses this computation to characterize which binomial edge ideals are Gorenstein \cite{gonzalez2021gorenstein}. The main goal of this paper is to investigate the computation of the log canonical threshold and F-threshold of a binomial edge ideal in terms of combinatorial properties of the underlying graph. To this end we introduce the following linear programming problem which is a modification of the linear programming problem considered by Shibuta and Takagi (see \cite{shibuta2009logcanonical}) and Hern\'{a}ndez (see \cite{hernandez2014Fpure}).

For $U \subset V(G)$, define $G[U]$ to be the induced graph on $U$. Let $\A := \{U \subset V(G) \mid \card{U} \geq 2\}$. Consider the following linear program associated to the graph $G$.
\begin{equation}
\begin{aligned}
\label{alg_path_packing_int_program}
\maximize \,  Z_{\PP,G} &:= \phantom{   a} \sum_{e\in E(G)} X_{\PP,e}  \\
\subject \phantom{lot of } & \\
  X_{\PP,e} & \in \k_{\geq 0} \text{ for all } e \in E(G) \\
   C_{\PP,U} &\leq \card{U} - 1  \text{ for all } U \in \A \\
   C_{\PP,v} &\leq 2  \text{ for all } v \in V(G)
\end{aligned}
\end{equation}
where for each $U \in \A$ and $v \in V$, $C_{\PP,U}$ and $C_{\PP,v}$ are the functions
\begin{align*}
 C_{\PP,U} &:= \phantom{ } \sum_{e \in E(G[U])} X_{\PP,e} \\
 C_{\PP,v} &:= \sum_{\substack{v \text{ is incident }\\ \text{with } e}} X_{\PP,e}.
\end{align*}
Throughout this paper let $\k$ denote a choice of ring $\Z$, $\Q$, or $\R$. Use \eqref{alg_path_packing_int_program} on different choices of $\k$, and denote the optimal value of \eqref{alg_path_packing_int_program} with respect to a choice of $\k$ by $\PP_{G,\k}$.

This paper studies $\PP_{G,\k}$ as it relates to the combinatorial properties of $G$ and the algebraic properties of $J_{G}$. In Section \ref{sec_Background} we introduce the relevant background. In Section \ref{sec_linear_progs_bin_edge_ideals} we show 
\begin{theoremx}
\label{thm_intro_inequalities_of_sec_two}
With the notation as above,
\begin{enumerate}
\item $\PP_{G,\Z}$ equals the binomial grade of $J_{G}$ (which is equal to the maximal number of edges in a semi-path of $G$)  (Proposition \ref{prop_PP_G,Z_counts_path_packings} and Remark \ref{rem_PP_G,Z_binomial_grade}),
\item The optimal value of the dual program over $\Z$, $\PP^{\ast}_{G,\Z}$, equals $\hgt J_{G} $ (Proposition \ref{prop_PP_dual_Z_equals_hgt}),
\item $\ft(J_{G,p}) \leq \PP_{G,\Q}$ for any prime number $p$ (Proposition \ref{prop_ft_leq_LP_G_Q}).
\end{enumerate}
\end{theoremx}
This establishes an analogue of Theorem \ref{thm_classical} parts \eqref{item_thm_classical_1} and \eqref{item_thm_classical_3}.

In Section \ref{sec_proof_cjec_block_graphs} we prove the main result of this paper.
\begin{theoremx}[Theorem \ref{thm_formula_ft_block_graph}, Corollary \ref{cor_main_conj_true_block_graphs}]
\label{thm_intro_ft_G_equals__PP_G,Q}
Let $G$ be a block graph, then
\begin{align*}
\ft(J_{G,p}) = \lct(J_{G}) = \PP_{G,\Q} \in \Z[\tfrac{1}{2}]
\end{align*}
for any prime number $p$.
\end{theoremx}
This result establishes a special case of an analogue of Theorem \ref{thm_classical} part \eqref{item_thm_classical_2}. 

Motivated by Theorem \ref{thm_intro_ft_G_equals__PP_G,Q}, Corollary \ref{prop_formula_ft_Konig_type}, and computational experiments we formulate
\bcjec
\label{conjectue_fpt_equal_LP}
Let $G$ be a graph, then
\begin{enumerate}
\item $\lct(J_{G}) = \PP_{G,\Q}$
\item  $\PP_{G,\Q} \in \Z[\frac{1}{2}]$.
\end{enumerate}
\ecjec
\noindent We prove this conjecture when $G$ is a block graph or $G$ is of K\"onig type. This conjecture suggests that a partial analogue of Theorem \ref{thm_classical} item \eqref{item_thm_classical_2} may hold.

In Section \ref{sec_ideals_of_Konig_type} we give a new characterization of graphs of K\"onig type, Lemma \ref{lem_characterization_Konig_type}. This enables us to show that every tree is of K\"onig type, Theorem \ref{thm_alg_constructing_T_correct}. This establishes an analogue of Theorem \ref{thm_konig}. In Section \ref{sec_ft_and_Ham_paths} we investigate for a connected graph $G$ the relationship between the equality $\PP_{G,\Q} = \card{V} - 1$ and $G$ being traceable.

In upcoming work \cite{laclair2023prep}, the author shows when the graph is a tree that the polytope determined by the constraints of Algorithm \ref{alg_path_packing_int_program} (i.e. fractional path packing polytope) has integral vertices, i.e. each coordinate of such a vertex is either $0$ or $1$. This result gives a partial analogue of Theorem \ref{thm_vertices_matching_poly}. The author also investigates the question whether for a general graph each coordinate of a vertex of the fractional path packing is either $0$, $1$, or $\frac{1}{2}$ \cite{laclair2023prep}.

\subsection*{Acknowledgments}
The author found computations in Macaulay 2 \cite{M2}, Python \cite{python}, Sage \cite{sage}, and Online Linear Optimization Solver \cite{onlineLPSolver} immensely helpful. The author would like to thank Alex Black, Jean-Philippe Labbe, and Vic Reiner for helpful conversations, Santiago Encinas for providing source code for the algorithm described in his paper \cite{blanco2018procedure}, and Uli Walther for many helpful discussions and suggestions.


\section{Background}

\label{sec_Background}

\subsection{Linear Programs}

For further background on linear programs and for proofs not contained here see \cite{dantzig2016linear}. A linear programming problem refers to the computation of the extremal value of a linear function over a convex polytope, i.e. any region defined by linear half-spaces. Any linear programming problem can be put into standard form which in this paper we take to be the following setup.

Let $A$ be an $n \times m$ matrix having entries in $\k$. Let $\mathbf{c} \in \k^{m}$ and $\mathbf{b} \in \k^{n}$. Then, the linear programming problem in standard form refers to
\begin{equation}
\label{eqn_primal}
\begin{split}
\maximize\limits_{\mathbf{x} \in \k^{m}} \quad & \mathbf{c}^{T} \mathbf{x} \\
\subject \hspace{.7em} &
  \phantom{ } \\
 & A \mathbf{x}  \leq \mathbf{b} \\
 & \phantom{A} \mathbf{x} \geq 0
\end{split}
\end{equation}
A linear programming problem written in this form will be called {\it primal}, and it can be shown that every linear programming problem is equivalent to a primal linear programming problem. 

Given a primal problem we can associate the following linear programming problem
\begin{equation}
\label{eqn_dual}
\begin{split}
\minimize_{\mathbf{y} \in \k^{n}} \quad & \mathbf{b}^{T} \mathbf{y} \\
\subject \hspace{.5em} &
  \phantom{ } \\
 & A^{T} \mathbf{y}  \geq \mathbf{c} \\
 & \phantom{A^{T}} \mathbf{y}  \geq 0
\end{split}
\end{equation}
which is called the {\it dual} linear programming problem of the primal.

Given a linear programming problem, a vector is called {\it feasible} if it satisfies the constraints of the problem. If such a vector exists the linear programming problem is called {\it feasible}. A linear programming problem is called {\it bounded} if the optimal value of the linear programming problem is finite. When a linear programming problem is bounded feasible we call the value attained by the linear program to be the {\it optimal value}, and we call a vector which realizes the optimal value and satisfies the constraints of the linear programming problem an {\it optimal solution}.

\bthm[{\cite[p.125 Duality Theorem]{dantzig2016linear}} Strong duality of linear programs]
\label{thm_strong_duality}

Given a primal (resp. dual) feasible, bounded linear programming problem, then the dual (resp. primal) programming problem is feasible, bounded and the optimal values of the primal (resp. dual) agrees with the optimal value of the dual (resp. primal) provided that $\Q \subset \k \subset \R$.
\ethm

\brem
When $\k = \Z$, then optimal value of a bounded feasible primal can be strictly smaller than the optimal value of the dual.
\erem

In the above setup of the Equation \eqref{eqn_primal}, let $A_{i}$ denote the $i$-th column of $A$ and $a_{j}^{T}$ denote the $j$-th row of $A$.

\bthm[ {\cite[p.136 Theorem 4]{dantzig2016linear}}  Complementary Slackness]
\label{thm_comp_slackness}
Let $\mathbf{x}$ and $\mathbf{y}$ be primal (respectively dual) optimal solutions of \eqref{eqn_primal} (respectively \eqref{eqn_dual}) and suppose that $\mathbf{c}^{T}\mathbf{x} = \mathbf{b}^{T} \mathbf{y}$, then
\begin{enumerate}
\item for $1 \leq i \leq m$, $(\mathbf{c}_{i} - \mathbf{y}^{T} A_{i})x_{i} = 0$. \label{item_comp_slackness_1}
\item for $1 \leq j \leq n$, $(a_{j}^{T} \mathbf{x} - \mathbf{b}_{j}) y_{j} = 0$. \label{item_comp_slackness_2}
\end{enumerate}
\ethm

\begin{pfsketch}
It is always true that
\begin{align*}
\mathbf{c}^{T} \mathbf{x} \leq \mathbf{y}^{T} A\mathbf{x} \leq \mathbf{y}^{T} \mathbf{b}.
\end{align*}
Equality of the outer two terms implies equality throughout and complementary slackness follows by distributing terms appropriately.
\end{pfsketch}

\subsection{Binomial Edge Ideals}

Throughout this paper we will only consider finite simple graphs, i.e. graphs without multiple edges or loops and on a finite number of vertices. Given a graph $G$ we denote its set of vertices by $V(G)$, or $V$ when the context is clear, and denote its set of edges by $E(G)$, or $E$ when the context is clear. Throughout this paper we will assume that $V(G) = \{1,\ldots,n\}$ for some positive integer $n$. Given a subset $U \subset V(G)$, we denote by $G[U]$ the subgraph of $G$ induced on the vertices $U$.

By a {\it path} in $G$ we mean an alternating sequence of distinct vertices and edges of $G$, $v_{1},e_{1},v_{2},\ldots,e_{\ell-1},v_{\ell}$ so that $e_{i} = \{v_{i},v_{i+1}\} \in E(G)$ for $1 \leq i \leq \ell-1$. By a {\it cycle} in $G$ we mean an alternating sequence of distinct vertices and edges of $G$, $v_{1},e_{1},v_{2},\ldots,e_{\ell-1},v_{\ell},e_{\ell}$ so that $e_{i} = \{v_{i},v_{i+1}\} \in E(G)$ for $1 \leq i \leq \ell-1$ and $e_{\ell} = \{v_{1},v_{\ell}\} \in E(G)$. By a {\it semi-path} we mean a subgraph $P$ of $G$ such that each connected component of $G$ is a path. Let $\hat{G}$ be the induced graph on $V(G) \setminus \mathcal{I}(G)$ where $\mathcal{I}(G)$ denotes the isolated vertices of $G$.

The following proposition gives a combinatorial description of the minimal primes of a binomial edge ideal and a formula for the height of prime ideals.

\bprop[{\cite[Lemma 3.1, Corollary 3.9]{herzog2010binomial}}]
\label{prop_height_min_primes_bin_edge_ideal}
Let $G$ be a graph on $n$ vertices. If $T \subset [n] := \{1,\ldots,n\}$, we will denote the induced subgraph on the vertices $G \setminus T$ by $G_{T}$. Let $c_{G}(T)$ denote the  number of connected components of $G_{T}$, and let $G_{1}, \ldots, G_{c_{G}(T)}$ denote the distinct connected components of $G_{T}$. Let $\tilde{G}_{i}$ denote the complete graph on the vertices $V(G_{i})$. Put 
\begin{align*}
P_{T}(G) := \big( \bigcup_{i\in T} \{X_{i},Y_{i}\}, J_{\tilde{G}_{1}}, \ldots, J_{\tilde{G}_{c_{G}(T)}} \big)
\end{align*}
where $J_{\tilde{G}_{\ast}}$ is as defined in Definition \ref{defn_bin_edge_ideal}.
Then, $P_{T}(G)$ is a minimal prime ideal of $J_{G}$ and 
\begin{align*}
\hgt P_{T}(G) = 2\card{T} + \sum_{i=1}^{c_{G}(T)} \left( \card{V(G_{i})} -1 \right) = n - c_{G}(T) + \card{T}.
\end{align*}
Moreover, every minimal prime ideal of $J_{G}$ is of the form $P_{T}(G)$ where $T = \varnothing$, or $T \neq \varnothing$ and for each $i \in T$ one has that $c_{G}(T\setminus i) < c_{G}(T)$.
\eprop

\bdefn
\label{defn_r(G)}
For a graph $G$ define the invariant \begin{align*}
r(G) := \sup \left\{ \card{E(P)} \mathrel{}\middle|\mathrel{} P \subset G \text{ is a semi-path} \right\}.
\end{align*}
We call a disjoint union of paths realizing $r(G)$ a {\it path packing of $G$}.
\edefn

\bdefn
\label{defn_traceable_graphs}
A graph $G$ is called {\it traceable} if $G$ has a Hamiltonian path.
\edefn

\brem
$G$ is traceable if and only if $r(G) = \card{V} - 1$.
\erem

\subsection{F-Thresholds}

\bdefn[{\cite[Lemma 1.1, Remark 1.5]{mustata2005FThresholds}}]
\label{defn_F-thresholds}
Let $S = \K[X_{1},\ldots, X_{n}]$ be a polynomial ring, $\fchar(\K) = p > 0$, $\m = (X_{1},\ldots, X_{n})$ the homogeneous maximal ideal, and $I$ an $S$-ideal with $I \subset \m$. For $e \in \N$, set
\begin{align*}
\nu_{e}(I) := \max\{r \mid I^{r} \not \subseteq \m^{[p^{e}]} \}
\end{align*}
where $\m^{[p^{e}]} := (X_{1}^{p^{e}},\ldots,X_{n}^{p^{e}})$. 
Then, 
\begin{align*}
\ft(I) := \lim_{e\ra \infty} \frac{\nu_{e}(I)}{p^{e}}
\end{align*}
exists and is called the {\it F-threshold} of $I$.
\edefn

We review three lemmas well-known to experts. The first lemma says that the computation of F-threshold is additive over ideals in disjoint variables and enables us to reduce the computation of the F-threshold of a binomial edge ideal to the computation of the F-threshold of the binomial edge ideal of its connected components. 

\blem
\label{lem_ft_additive}
Let $\fchar \K = p > 0$, $I \subset S := \K[X_{1},\ldots,X_{n}]$ with $\mathfrak{m}$ its homogeneous maximal ideal, $J \subset S^{'} := \K[Y_{1},\ldots,Y_{m}]$ with $\mathfrak{n}$ its homogeneous maximal ideal, $T := S \tensor_{\K} S^{'}$ with $\mathfrak{o}$ its homogeneous maximal ideal. Then, 
\begin{align*}
\ft_{T}(I\cdot T + J \cdot T) = \ft_{S}(I) + \ft_{S^{'}}(J).
\end{align*}
In particular, if $G$ is a graph and $G_{1},\ldots, G_{c}$ are its connected components, then 
\begin{align*}
\ft(J_{G}) = \sum_{i=1}^{c} \ft(J_{G_{i}}).
\end{align*}
\elem

The next lemma says that the F-threshold is monotonic along ideal inclusion which allows us to obtain lower bounds on a graph via its subgraphs.

\blem
\label{lem_ft_monotonic}
Let $I \subset J \subset \mathfrak{m}$ be ideals in $\K[X_{1},\ldots,X_{n}]$. Then, 
\begin{align*}
\ft(I) \leq \ft(J).
\end{align*} 
In particular, if $H$ is a subgraph of $G$, then 
\begin{align*}
\ft(J_{H}) \leq \ft(J_{G}).
\end{align*}
\elem

\bproof
It follows from the definition that $\nu_{e}(I) \leq \nu_{e}(J)$ for all $e \in \N$.
\eproof

The final lemma says that the F-threshold can only decrease under taking initial ideals.

\blem[{\cite[Proposition 4.5]{takagi2004fpure}}]
\label{lem_ft_decr_taking_initial_ideal}
Let $I \subset \mathfrak{m}$ be ideal in $\K[X_{1},\ldots,X_{n}]$ and $<$ a term order on the polynomial ring. Then, $\ft(\init_{<}(I)) \leq \ft(I)$.
\elem

Miller, Singh, and Varbaro \cite{miller2014fptDeterminantal} computed the F-threshold of a determinantal ideal of a generic matrix as follows.

\begin{thrm}[{\cite[Theorem 1.2]{miller2014fptDeterminantal}}]
\label{thm_msv_ft_determinantal_ideal}
Fix positive integers $t \leq m \leq n$, and let $X$ be an $m \times n$ matrix of indeterminates over a field $\mathbb{F}$ of prime characteristic. Let $S$ be the polynomial ring $\mathbb{F}[X]$, and $I_{t}$ the ideal generated by the size $t$ minors of $X$.

The F-threshold of $I_{t}$ is
\begin{align*}
\min\bigg\lbrace \frac{(m-k)(n-k)}{t-k} \mid k = 0,\ldots, t-1 \bigg\rbrace.
\end{align*}
\end{thrm}

\brem
\label{rem_ft_of_complete_graph}
In the above theorem, when $t = m = 2 \leq n$, then $I_{t}(X)$ corresponds to $J_{K_{n}}$ where $K_{n}$ denotes the complete graph on $n$ vertices, and thus 
\begin{align*}
\ft(J_{K_{n}}) = n-1.
\end{align*}
\erem

In the master's thesis of Velez \cite{velez2019f}, he proves 
\bprop[{{\cite[Theorem 3.1.6]{velez2019f}}}]
\label{prop_velez_computation}
Let $M$ be a generic $m \times n$ matrix with $m \leq n$. Let $\delta_{i}$ denote the determinant of the submatrix given by the columns $\{i,i+1,\ldots,i+m-1\}$ for $1 \leq i \leq n-m+1$. Let $J := ( \delta_{i} \mid 1 \leq i \leq n-m+1)$. Then, $\ft(J) = n-m+1$.
\eprop

\begin{pfsketch}
By Lemma \ref{lem_ft_decr_taking_initial_ideal} and Theorem \ref{thm_msv_ft_determinantal_ideal},
\begin{align*}
n-m+1 = \ft(\init_{<}(J_{G})) \leq \ft(J_{G}) \leq n-m+1.
\end{align*}
\end{pfsketch}

\brem
\label{rem_ft_path}
When $t = m = 2$ in Proposition \ref{prop_velez_computation} the ideal $J_{G}$ corresponds the binomial edge ideal of a path on $n$ vertices. I.e., if $G = P_{n}$ is a path on $n$ vertices, then $\ft(J_{G}) = n-1$.
\erem

\bthm[\cite{lucas1878theorie} Lucas's theorem]
\label{thm_Lucas}
Let $p$ be a prime integer and $m$ and $n$ positive integers. Choose a positive integer $e$ so that we can write $m = \sum_{i=0}^{e} m_{i}p^{i}$ and $n = \sum_{i=0}^{e} n_{i} p^{i}$ where $0 \leq m_{i}, n_{i} < p$ for $0 \leq i \leq e$ . Then,
\begin{align*}
\binom{n}{m} \equiv \prod_{i=1}^{e} \binom{n_{i}}{m_{i}} \mod{p}
\end{align*}
with the convention that $\binom{a}{b} = 0$ if $a < b$.
\ethm

\section{Linear Programs and Binomial Edge Ideals}

\label{sec_linear_progs_bin_edge_ideals}

\subsection{Binomial Grade, Height, and Linear Programs}

Recall that $r(G)$ is the maximal length of a path packing (Definition \ref{defn_r(G)}), and $\PP_{G,\Z}$ is the optimal value of an integer program (discussion surrounding \eqref{alg_path_packing_int_program}).

\bprop
\label{prop_PP_G,Z_counts_path_packings}
Given a graph $G$, $r(G) = \PP_{G,\Z}$.
\eprop

\bproof
For $v \in V$, the constraint $C_{\PP,v} \leq 2$ forces that a feasible solution has at most two edges incident to $v$. For $U \in \A$, the constraint $C_{\PP,U} \leq \card{U}-1$ forces that a feasible solution does not induce a Hamiltonian cycle on $G[U]$. Thus, an optimal solution to \eqref{alg_path_packing_int_program} consists of disjoint paths. Conversely, a path packing is a feasible solution of \eqref{alg_path_packing_int_program}.
\eproof

\brem
\label{rem_PP_G,Z_binomial_grade}
The feasible solutions of \eqref{alg_path_packing_int_program} are precisely the semi-paths of $G$. Moreover, in \cite[Lemma 3.2, Lemma 3.3]{herzog2022graded} the authors show that $r(G)$ is equal to the binomial grade of $J_{G}$, i.e. the length of a longest regular sequence in $J_{G}$ consisting of binomials which form part of a minimal generating set of $J_{G}$.
\erem

The dual linear program of \eqref{alg_path_packing_int_program} solved over $\k$ is given by
\begin{equation}
\label{alg_dual_path_packing_int_program}
\begin{split}
\minimize \quad & Z_{\PPd,G} := \sum_{v \in V(G)} 2 \cdot W_{\PPd,v} + \sum_{U \in \A} (\card{U} - 1) \cdot W_{\PPd,U} \\
\subject  \hspace{.5em} &  \phantom{ } \\
 & W_{\PPd,v} \in \k_{\geq 0} \text{ for all } v \in V(G) \\
 & W_{\PPd,U} \in \k_{\geq 0} \text{ for all }  U \in \A \\
 & D_{\PPd,e} := W_{\PPd,i} + W_{\PPd,j} + \sum_{\substack{U \in \A  \\ e \in G[U] } } W_{\PPd,U} \geq 1 \text{ for all } e = \{i,j\} \in E(G)
\end{split}
\end{equation}
Denote the optimal value of \eqref{alg_dual_path_packing_int_program} over $\k$ by $\PP_{G,\k}^{\ast}$.

From a combinatorial perspective, the integral linear program \eqref{alg_dual_path_packing_int_program} computes a minimal \textit{weighted} edge covering of the graph by stars and induced subgraphs. From a commutative algebra perspective, \eqref{alg_dual_path_packing_int_program} computes the height of $J_{G}$.

\blem
\label{lem_nice_form_optimal_soln_dual_path_packing}
Fix a graph $G$ and set $\k = \Z$. There exists an optimal solution of the linear program \eqref{alg_dual_path_packing_int_program} of the form $( \{b_{v} \}_{v\in V}, \{c_{U}\}_{U\in\A} )$ satisfying
\begin{enumerate}
\item $b_{v} \in \{0,1\}$ for all $v \in V$ and $c_{U} \in \{0,1\}$ for all $U \in \A$,
\item if $U$ and $U^{'}$ are elements of $\mathcal{A}$ with $c_{U} = c_{U^{'}} = 1$, then $U \cap U^{'} = \varnothing$,
\item if $b_{v} = 1$ for some $v \in V(G)$, then $v \notin U$ for all $c_{U} =1$,
\item if $c_{U}  = 1$ for some $U \in \A$, then $U$ is connected.
\end{enumerate}
For such a solution, let $T:= \{v \in V \mid c_{v} = 1\}$ and $U_{1},\ldots,U_{c}$ to be the elements of $\A$ which satisfy $c_{U_{i}} = 1$. Then, the $U_{i}$ are the connected components of $\wh{G \setminus T}$.
\elem

\bproof
Let $( \{b_{v} \}_{v\in V}, \{c_{U}\}_{U\in\A} )$ be any optimal solution to \eqref{alg_dual_path_packing_int_program}.
\begin{enumerate}
\item We may assume that $b_{v}, c_{U} \in \{0,1\}$ for all $v \in V$ and $U \in \A$. Otherwise, for each non-zero $b_{v}$ or non-zero $c_{U}$ reassign it the value $1$. 
\item We may suppose that if $U$ and $U^{'}$ are elements of $\A$ with $c_{U} = c_{U^{'}} = 1$, then $U \cap U^{'} = \varnothing$. Otherwise, redefine $c_{U} = c_{U^{'}} = 0$ and $c_{U\cup U^{'}} = 1$.
\item We may suppose that if $b_{v} = 1$ for some $v \in V(G)$, then $v \notin U$ for all $c_{U} =1$. Otherwise, if $\card{U} > 2$ put $c_{U} = 0$ and $c_{U \setminus \{v\}} = 1$, and if $\card{U} = 2$ put $c_{U} = 0$. 
\item We may suppose that if $c_{U}  = 1$ for some $U \in \A$, then $U$ is connected. Otherwise, write $U = U_{1} \cup U_{2}$ with $U_{1} \cup U_{2} = \varnothing$. Put $c_{U} = 0$ and $c_{U_{i}} =1$ if $U_{i} \in \A$ for $i = 1,2$.
\end{enumerate}
After each application of one of the above steps this defines a new feasible solution of \eqref{alg_dual_path_packing_int_program} which does not increase the cost function.

After these reductions, define $T := \{v \in V \mid a_{v} = 1\}$ and $U_{1},\ldots,U_{c}$ to be elements of $\A$ which satisfy $c_{U_{i}} = 1$. The above reductions together with the constraints $D_{\PPd,e} \geq 1$ for every $e \in E(G)$ implies that the $U_{i}$ are the connected components of $\wh{G \setminus T}$.
\eproof

\bprop
\label{prop_PP_dual_Z_equals_hgt}
Let $G$ be a graph, then $\PP_{G,\Z}^{\ast} = \hgt J_{G}$.
\eprop

\bproof
We first show $\PP_{G,\Z}^{\ast} \leq \hgt J_{G}$. Let $T \subset V(G)$ be chosen so that $\hgt J_{G} = \hgt P_{T}(G) = 2 \card{T} + \sum_{i=1}^{c_{T}(G)} \big( \card{V(G_{i})} - 1 \big)$ where the notation is as in Proposition \ref{prop_height_min_primes_bin_edge_ideal}. For $v \in V$ let
\begin{align*}
b_{v} := 
\begin{cases}
 1, & v \in T \\
 0, & \text{else}
\end{cases}
\end{align*}
and for $U \in \A$ let
\begin{align*}
c_{U} = 
\begin{cases}
 1, & U = V(G_{i}) \text{ for some } 1 \leq i \leq c_{T}(G)\\
 0, & \text{else}
\end{cases}.
\end{align*}
Then the tuple $\big( \{b_{v}\}_{v \in V}, \{c_{U}\}_{U \in \A} \big)$ is a feasible solution to \eqref{alg_dual_path_packing_int_program}, and moreover
\begin{align*}
\PP_{G,\Z}^{\ast} \leq Z_{\PPd,G} \big( \{b_{v}\}_{v \in V}, \{c_{U}\}_{U \in \A} \big) = 2 \card{T} + \sum_{i=1}^{c_{T}(G)} \big( \card{V(G_{i})} - 1 \big) = \hgt J_{G}.
\end{align*}
Next, we show that $\hgt J_{G} \leq \PP_{G,\Z}^{\ast}$. Let $\big( \{b_{v}\}_{v \in V}, \{c_{U}\}_{U \in \A} \big)$ be an optimal solution of \eqref{alg_dual_path_packing_int_program} when $\k = \Z$ of the form guaranteed by Lemma \ref{lem_nice_form_optimal_soln_dual_path_packing}. Let $T$ and $U_{1},\ldots, U_{c}$ be as constructed by Lemma \ref{lem_nice_form_optimal_soln_dual_path_packing}. The cost of this optimal tuple is precisely $2 \card{T} + \sum_{i=1}^{c} (\card{U_{i}} -1)$ which is equal to $\hgt P_{T}(G)$. This proves that 
\begin{align*}
\hgt J_{G} \leq \hgt P_{T}(G) = \PP_{G,\Z}^{\ast}
\end{align*}
which completes the proof.
\eproof

\subsection{Thresholds and Linear Programs}

Next, we study the LP relaxation of the above linear programs, i.e. we allow the variables to take values in $\R_{\geq 0}$, and how the optimal solution relates to commutative algebra invariants. 

For a graph $G$ denote by $\LP_{G,\k}$ the optimal value of the following linear program over $(\k_{\geq 0})^{2\card{E}}$.

\begin{equation}
\label{alg_primal_LP_G}
\begin{split}
\maximize \quad & Z_{\LP,G} := \sum_{e\in E(G)} (X_{\LP,e} + Y_{\LP,e}) \\
\subject \hspace{.7em} &  \phantom{ } \\
 & X_{\LP,e} \in \k_{\geq 0} \text{ for all } e \in E(G) \\
  & Y_{\LP,e} \in \k_{\geq 0} \text{ for all } e \in E(G) \\
  & C_{\LP,U} := \sum_{e \in E(G[U])} (X_{\LP,e} + Y_{\LP,e}) \leq \card{U} - 1  \text{ for all } U \in \A \\
  & C_{\LP,v,x} := \sum_{\substack{e = \{v,w\}\\ v < w}  } X_{\LP,e} + \sum_{\substack{e = \{w,v\}\\ w < v}  } Y_{\LP,e} \leq 1  \text{ for all } v \in V(G) \\
  &  C_{\LP,v,y} := \sum_{\substack{e = \{v,w\}\\ v < w}  } Y_{\LP,e} + \sum_{\substack{e = \{w,v\}\\ w < v}  } X_{\LP,e} \leq 1 \text{ for all } v \in V(G) \\  
\end{split}
\end{equation}

\brem
\label{rem_LP_G,Q_reduces_connected_cases}
Let $G$ be a graph and $G_{1}, \ldots, G_{c}$ the connected components of $G$. Then,
\begin{align*}
\LP_{G,\Q} = \sum_{i=1}^{c} \LP_{G_{i},\Q}.
\end{align*}
\erem

The following proposition says that the constraint $C_{\LP,U} \leq \card{U} - 1$ in \eqref{alg_primal_LP_G} can be strengthened by requiring that in addition $U$ is biconnected. Recall that a connected graph $G$ is biconnected if $G\setminus \{v\}$ is connected for every $v \in V(G)$. In the following linear program we differentiate the intdeterminates and constraints from those appearing in \eqref{alg_primal_LP_G} by utilizing the superscript ``$\bi$".

\bprop
The set of feasible solutions of \eqref{alg_primal_LP_G} is the same as the set of feasible solutions of the following linear program
\begin{equation}
\label{alg_primal_LP_G_biconnected}
\begin{split}
\maximize \quad & Z_{\LP,G}^{\bi} := \sum_{e\in E(G)} (X_{\LP,e}^{\bi} + Y_{\LP,e}^{\bi}) \\
\subject \hspace{.7em} &  \phantom{ } \\
 & X_{\LP,e}^{\bi} \in \k_{\geq 0} \text{ for all } e \in E(G) \\
  & Y_{\LP,e}^{\bi} \in \k_{\geq 0} \text{ for all } e \in E(G) \\
  & C_{\LP,U}^{\bi} := \sum_{e \in E(G[U])} (X_{\LP,e}^{\bi} + Y_{\LP,e}^{\bi}) \leq \card{U} - 1  \text{ for all } U \in \A, \phantom{ } G[U] \text{ biconnected} \\
  & C_{\LP,v,x}^{\bi} := \sum_{\substack{e = \{v,w\}\\ v < w}  } X_{\LP,e}^{\bi} + \sum_{\substack{e = \{w,v\}\\ w < v}  } Y_{\LP,e}^{\bi} \leq 1  \text{ for all } v \in V(G) \\
  &  C_{\LP,v,y}^{\bi} := \sum_{\substack{e = \{v,w\}\\ v < w}  } Y_{\LP,e}^{\bi} + \sum_{\substack{e = \{w,v\}\\ w < v}  } X_{\LP,e}^{\bi} \leq 1 \text{ for all } v \in V(G) \\  
\end{split}
\end{equation}
where $v < w$ means that the label on the vertex $v$ is smaller than the label on the vertex $w$.
\eprop

\bproof
It suffices to show that feasible solutions of \eqref{alg_primal_LP_G_biconnected} are feasible solutions of \eqref{alg_primal_LP_G}. It suffices to restrict to the case of $U \in \A$ is connected. Write $U = U_{1} \cup \cdots \cup U_{c}$ as a union of its maximal biconnected components. Let $\mbf{a}$ denote a feasible solution of \eqref{alg_primal_LP_G}. Observe that 
\begin{align*}
C_{U}(\mbf{a}) = \sum_{i=1}^{c} C_{U_{i}}(\mbf{a}) = \sum_{i=1}^{c} C_{U_{i}}^{\bi}(\mbf{a}) \leq \sum_{i=1}^{c} (\card{U_{i}} -1).
\end{align*}
The above assumptions on $U$ together with induction on $c$ shows that 
\begin{align*}
\sum_{i=1}^{c} (\card{U_{i}} -1) \leq \card{U} -1.
\end{align*}
\eproof

\brem
Linear program \eqref{alg_primal_LP_G_biconnected} will have $O(2^{n})$ constraints for a complete graph on $n$ vertices.
\erem

\bex
Let $G$ be the graph depicted in Figure \ref{fig_the_net}.

\def\putBelow{-.5*.9}
\def\edgeDist{1.75*.9}
\def\rtTriangleEdgeMultiplier{1/8*.9}

\begin{figure}[h]
\begin{center}
\begin{tikzpicture}
\filldraw[black] (-3*\edgeDist,0) circle (2pt) node at (-3*\edgeDist,\putBelow) {$1$};
\filldraw[black] (-2*\edgeDist,0) circle (2pt) node at (-2*\edgeDist,\putBelow) {$2$};
\filldraw[black] (-1*\edgeDist,0) circle (2pt) node at (-1*\edgeDist,\putBelow) {$4$};
\filldraw[black] (0,0) circle (2pt) node at (0,\putBelow) {$6$};
\filldraw[black] (-3/2*\edgeDist, 
0.866025403784*\edgeDist) circle (2pt) node at (-3/2*\edgeDist+.5, 
0.866025403784*\edgeDist) {$3$};
\filldraw[black] (-3/2*\edgeDist,1.75*\edgeDist) circle (2pt) node at (-3/2*\edgeDist+.5,1.75*\edgeDist) {$5$};

\draw (-3*\edgeDist,0) to (-2*\edgeDist,0);
\draw (-2*\edgeDist,0) to (-1*\edgeDist,0);
\draw (-1*\edgeDist,0) to (0,0);
\draw (-2*\edgeDist,0) to (-3/2*\edgeDist, 
0.866025403784*\edgeDist);
\draw (-1*\edgeDist,0) to (-3/2*\edgeDist, 
0.866025403784*\edgeDist);
\draw (-3/2*\edgeDist, 
0.866025403784*\edgeDist) to (-3/2*\edgeDist,1.75*\edgeDist);

\node at (-5/2*\edgeDist,-.25) {a};
\node at (-3/2*\edgeDist,-.25) { c};
\node at (-1/2*\edgeDist,-.25) { f};
\node at (-3/2*\edgeDist+.25,1.375*\edgeDist) { e};

\node at (-1.25*\edgeDist +
1.73205080757*\rtTriangleEdgeMultiplier, 0.866025403784/2*\edgeDist + \rtTriangleEdgeMultiplier) { d};

\node at (-1.75*\edgeDist -
1.73205080757*\rtTriangleEdgeMultiplier, 0.866025403784/2*\edgeDist + \rtTriangleEdgeMultiplier) { b};

\end{tikzpicture}
\end{center}
\caption{Net}
\label{fig_the_net}
\end{figure}
We label the edges of the graph with the letters a through f and the vertices with the numbers $1$ through $6$ as in Figure \ref{fig_the_net}. Then, the constraints in \eqref{alg_primal_LP_G_biconnected} are 
\begin{align*}
C_{\LP,\{1,2\}}^{\bi} &= X_{a}^{\bi} + Y_{a}^{\bi} \leq 1 & C_{\LP,4,x}^{\bi} &= X_{f}^{\bi} + Y_{c}^{\bi} + Y_{d}^{\bi} \leq 1 & C_{\LP,1,x}^{\bi} &= X_{a}^{\bi} \leq 1 \\
C_{\LP,\{2,3\}}^{\bi} &= X_{b}^{\bi} + Y_{b}^{\bi} \leq 1 & C_{\LP,4,y}^{\bi} &= Y_{f}^{\bi} + X_{c}^{\bi} + X_{d}^{\bi} \leq 1 & C_{\LP,1,y}^{\bi} &= Y_{a}^{\bi} \leq 1 \\
C_{\LP,\{2,4\}}^{\bi} &= X_{c}^{\bi} + Y_{c}^{\bi} \leq 1 & C_{\LP,2,x}^{\bi} &= Y_{a}^{\bi} + X_{b}^{\bi} + X_{c}^{\bi} \leq 1 & C_{\LP,5,x}^{\bi} &= Y_{e}^{\bi} \leq 1 \\
C_{\LP,\{3,4\}}^{\bi} &= X_{d}^{\bi} + Y_{d}^{\bi} \leq 1 & C_{\LP,2,y}^{\bi} &= X_{a}^{\bi} + Y_{b}^{\bi} + Y_{c}^{\bi} \leq 1 & C_{\LP,5,y}^{\bi} &= X_{e}^{\bi} \leq 1 \\
C_{\LP,\{3,5\}}^{\bi} &= X_{e}^{\bi} + Y_{e}^{\bi} \leq 1 & C_{\LP,3,x}^{\bi} &= X_{d}^{\bi} + X_{e}^{\bi} + Y_{b}^{\bi} \leq 1 & C_{\LP,6,x}^{\bi} &= Y_{f}^{\bi} \leq 1 \\
C_{\LP,\{4,6\}}^{\bi} &= X_{f}^{\bi} + Y_{f}^{\bi} \leq 1 & C_{\LP,3,y}^{\bi} &= Y_{d}^{\bi} + Y_{e}^{\bi} + X_{b}^{\bi} \leq 1 & C_{\LP,6,y}^{\bi} &= X_{f}^{\bi} \leq 1
\end{align*}
%
%
%
%
%
%
and $C_{\{2,3,4\}}^{\bi} = X_{b}^{\bi} + Y_{b}^{\bi} + X_{c}^{\bi} + Y_{c}^{\bi} + X_{d}^{\bi} + Y_{d}^{\bi} \leq 2$.
\eex

In the linear program \eqref{alg_dual_LP_G} below, whenever we write $\{i,j\}$ for an edge of a graph we will in addition assume that $i < j$. The dual of the linear program \eqref{alg_primal_LP_G} is given by 
\begingroup
\allowdisplaybreaks
\begin{equation}
\label{alg_dual_LP_G}
\begin{split}
\minimize \hspace{.85em} & Z_{\LPd,G} := \sum_{v \in V(G)} (W_{\LPd,v,x} + W_{\LPd,v,y}) + \sum_{U \in \A} (\card{U} -1) W_{\LPd,U} \\
\subject \hspace{.25em} &  \phantom{ } \\
& W_{\LPd,v,x} \in \k_{\geq 0} \text{ for all } v \in V(G) \\
& W_{\LPd,v,y} \in \k_{\geq 0} \text{ for all } v \in V(G) \\
& W_{\LPd,U} \in \k_{\geq 0} \text{ for all } U \subset V(G) \text{ with } \card{U} \geq 2 \\
& D_{\LPd,e,x} := \sum_{\substack{e \in E(G[U])\\U\in \A}} W_{\LPd,U} + W_{\LPd,i,x} + W_{\LPd,j,y} \geq 1 \text{ for all } e = \{i,j\} \in E(G) \\
& D_{\LPd,e,y} := \sum_{\substack{e \in E(G[U])\\U\in \A}} W_{\LPd,U} + W_{\LPd,j,x} + W_{\LPd,i,y} \geq 1 \text{ for all } e = \{i,j\} \in E(G)
\end{split}
\end{equation}
\endgroup
We denote by $\LP_{G,\k}^{\ast}$ the optimal value of \eqref{alg_dual_LP_G}.

\bprop
Let $G$ be a labeled graph. Then, $\PP_{G,\k} = \LP_{G,\k}$ and $\PP_{G,\k}^{\ast} = \LP_{G,\k}^{\ast}$.
\eprop

\bproof
Let $\Sigma$ (respectively $\Sigma^{'}$) denote the set of feasible solutions of \eqref{alg_path_packing_int_program} (respectively of the feasible solutions of \eqref{alg_primal_LP_G}). If $\k$ is a field, then there are linear bijections 
\begin{align*}
\Sigma &\lra \Sigma^{'} \\
(a_{e})_{e \in E(G)} &\bra \big( \{\frac{a_{e}}{2} \}_{e \in E(G)}, \{\frac{a_{e}}{2})\}_{e \in E(G)} \big) \\
(a_{e} + b_{e})_{e\in E(G)} &\mapsfrom \big( \{ a_{e} \}_{e \in E(G)}, \{ b_{e} \}_{e \in E(G)} \big)
\end{align*}
If $\k = \Z$, then there are a linear bijections
\begin{align*}
\Sigma &\lra \Sigma^{'} \\
(a_{e})_{e \in E(G)} &\bra \big( \{ a_{e} \}_{e \in E(G)}, \{0\}_{e \in E(G)} \big) \\
(a_{e} + b_{e})_{e\in E(G)} &\mapsfrom \big( a_{e} \}_{e \in E(G)}, \{ b_{e} \}_{e \in E(G)} \big)
\end{align*}
The constraint $C_{U} \leq 1$ for $U = e$ shows that the map from $\Sigma^{'} \ra \Sigma$ is well-defined.

The dual case is proved analogously.
\eproof

\bprop
\label{r(G)_lt_fpt}
Let $G$ be a graph, then
\begin{align*}
r(G) \leq \ft(J_{G})
\end{align*}
where $r(G)$ is as found in Definition \ref{defn_r(G)}.
\eprop

\bproof
Let $<$ denote the lex monomial order on $R$ with $X_{1} > \cdots > X_{n} > Y_{1} > \cdots > Y_{n}$. Then,
\begin{align*}
r(G) \leq \ft (\init_{<} J_{G}) \leq \ft (J_{G}).
\end{align*}
\eproof

\blem[{\cite[Lemma 1.2]{miller2014fptDeterminantal}}]
\label{lem_MSV_1.2}
Fix a $2 \times n$ generic matrix $M$ and $S = \mathbb{F}_{p}[M]$. There exists a non-negative integer $N_{M}$, so that for every $q = p^{e}$ where $e \in \Z_{\geq 0}$, we have $\nu_{q}(I_{2}(M)) \leq N_{M} + (q-1)(n-1)$.
\elem

\bproof
Follows from the proof of their Lemma 1.2 by taking, in their notation, $m = t = 2$ and $k = 1$.
\eproof

\bprop
\label{prop_ft_leq_LP_G_Q}
Let $G$ be a labeled graph, then 
\begin{align*}
\ft(J_{G}) \leq \LP_{G,\Q}
\end{align*}
\eprop

\bproof
Fix notation of $M$ as a $2 \times \card{V(G)}$ generic matrix, $S = \mathbb{F}_{p}[M]$. For $U \subset V(G)$ we consider $M_{U}$ and $S_{U}$ the restriction of $M$ and $S$ to the submatrix (respectively subring) of variables corresponding to the vertices of $U$. We take $N := \max\big\lbrace  \card{V(G)}, \{ N_{M_{U}} \mid U \in \A \} \big\rbrace$ where $N_{M_{U}}$ is as defined in Lemma \ref{lem_MSV_1.2}.

We denote by $P_{\infty}$ the polytope contained in $(\R_{\geq 0})^{2\card{E}}$ determined by the constraints of linear program \eqref{alg_primal_LP_G}.

For every $t \in \N$, let $P_{t}$ be the polytope inside $(\R_{\geq 0})^{2\card{E}}$ determined by the following constraints
\begin{align}
 & \tilde{X}_{\LP,e} \geq 0 \text{ for all } e \in E(G) \nonumber \\
  & \tilde{Y}_{\LP,e} \geq 0 \text{ for all } e \in E(G) \nonumber \\
  & \tilde{C}_{\LP,U} := \sum_{e \in E(G[U])} (\tilde{X}_{\LP,e} + \tilde{Y}_{\LP,e}) \leq \frac{N + (p^{t} -1)(\card{U} -1)}{p^{t}} \text{ for all } U \in \A \label{eqn_P_t_first} \\
  & \tilde{C}_{\LP,v,x} := \sum_{\substack{e = \{v,w\}\\ v < w}  } \tilde{X}_{\LP,e} + \sum_{\substack{e = \{w,v\}\\ w < v}  } \tilde{Y}_{\LP,e} \leq 1  \text{ for all } v \in V(G) \label{eqn_P_t_sec} \\
  &  \tilde{C}_{\LP,v,y} := \sum_{\substack{e = \{v,w\} \\ v < w}  } \tilde{Y}_{\LP,e} + \sum_{\substack{e = \{w,v\} \\ w < v}  } \tilde{X}_{\LP,e} \leq 1 \text{ for all } v \in V(G) \label{eqn_P_t_third}
\end{align}
We denote by $\LP_{G,\Q,t}$ to be the maximum of the function 
\begin{align*}
\R^{2\card{E}} \ra \R : (\{ \tilde{X}_{\LP,e} \}_{e \in E(G)}, \{ \tilde{Y}_{\LP,e} \}_{e \in E(G)} ) \bra \sum_{e \in E(G)} ( \tilde{X}_{\LP,e} + \tilde{Y}_{\LP,e})
\end{align*}
when restricted to the region $P_{t}$.

Fix $t$ a positive integer. Then, there exist non-negative integers $a_{e}, b_{e}$ for $e \in E(G)$ satisfying $\sum_{e \in E(G)} (a_{e} + b_{e}) = \nu_{t}(J_{G})$ and a monomial 
\begin{align*}
\mu := \prod\limits_{\substack{e = (i,j) \in E(G)\\ i < j}} (X_{i}Y_{j})^{a_{e}} (X_{j}Y_{i})^{b_{e}}
\end{align*}
 belonging to the monomial support of some element of $J_{G}^{\nu_{t}(J_{G})} \setminus m^{[p^{t}]}$.

Observe that the tuple $\left( \mbf{a}, \mbf{b} \right) = \left( \{\frac{a_{e}}{p^{t}} \}_{e \in E(G)}, \{ \frac{b_{e}}{p^{t}} \}_{e \in E(G)} \right) \in P_{t}$. Indeed, since $\mu \notin m^{[p^{t}]}$,
\begin{align*}
 \tilde{C}_{\LP,v,x}\big( \mbf{a}, \mbf{b} \big) &= \deg_{X_{v}} \mu \leq p^{t} \\
\tilde{C}_{\LP,v,y}\big( \mbf{a}, \mbf{b} \big) &= \deg_{Y_{v}} \mu \leq p^{t}.
\end{align*}
Indeed, by Lemma \ref{lem_MSV_1.2},
\begin{align*}
\tilde{C}_{\LP,U}\big( \mbf{a}, \mbf{b} \big) \leq N + (p^{t}-1)(\card{U} -1).
\end{align*}
This shows that $\frac{\nu_{t}(J_{G})}{p^{t}} \leq \LP_{G,\Q,t}$. Taking limit as $t \ra \infty$ it follows by Lemma \ref{lem_aux_lem} that $\ft(J_{G}) \leq \LP_{G,\Q}$.
\eproof

\blem
\label{lem_aux_lem}
With the notation as above,
\begin{enumerate}
\item For  every $0 < \epsilon < 1$, there exists $t \in \N$ so that for every $\mbf{x} \in P_{t}$ there exists $\mbf{y} \in P_{\infty}$ with $\norm{\mbf{x} - \mbf{y}}_{1} < \epsilon$. \label{item_aux_lem_1}
\item $\lim_{t \ra \infty} \LP_{G,\Q,t} = \LP_{G,\Q}$.
\end{enumerate}
\elem

\bproof
\begin{enumerate}
\item Fix $0 < \epsilon < 1$. Take $t \in \N$ so that $2 \card{E}\frac{N}{p^{t}} < \epsilon$. For $\mbf{x} \in P_{t}$, define $\mbf{y}$ via $y_{i} := x_{i} - \min\{x_{i},\frac{N}{p^{t}}\}$ for $1 \leq i \leq 2 \card{E}$. We observe that 
\begin{align*}
\norm{ \mbf{x} - \mbf{y} }_{1} \leq 2 \card{E} \frac{N}{p^{t}} < \epsilon.
\end{align*}
Since $\mbf{y}$ is pointwise less than $\mbf{x}$ and $\mbf{x}$ satisfies the constraints \eqref{eqn_P_t_sec} and \eqref{eqn_P_t_third}, $C_{\LP,v,x}(\mbf{y}) < 1$ and $C_{\LP,v,y}(\mbf{y}) < 1$ for all $v \in V$. Thus, to show that $\mbf{y} \in P_{\infty}$ it suffices to show that $C_{U}(\mbf{y}) \leq \card{U} -1$ for all $U \in \A$.

\bobs
\label{obs_one}
The constraint $\tilde{C}_{\LP,U}
(\mbf{x}) \leq \frac{N + (p^{t}-1)(\card{U}-1)}{p^{t}}$ implies that 
\begin{align}
\label{eqn_obs_in_proof_1}
\tilde{C}_{\LP,U}(\mbf{x}) - \frac{N}{p^{t}} < \card{U}-1.
\end{align}
\eobs

\bobs
\label{obs_two}
If $\tilde{C}_{\LP,U}(\mbf{x}) > \card{U} -1$, then $x_{\max,U} \geq \frac{N}{p^{t}}$ where  $x_{\max,U}$ is the largest value of $\mbf{x}$ that appears in the support of $\tilde{C}_{\LP,U}(\mbf{x})$. If this is not the case, then we have that
\begin{align*}
\card{U}-1 &< \tilde{C}_{\LP,U}(\mbf{x}) \\
&\leq 2 \card{E(G[U])} \cdot x_{\max,U} \\
&\leq 2 \card{E} \frac{N}{p^{t}} \\
& < \epsilon \\
& < 1.
\end{align*}
But, then this implies that $\card{U} < 2$; a contradiction.
\eobs

There is nothing to prove if $\tilde{C}_{\LP,U}(\mbf{x}) \leq \card{U} -1$. If $\tilde{C}_{\LP,U} > \card{U}-1$, then by Observation \ref{obs_two} there exists $x_{i}$ appearing in the support of $\tilde{C}_{\LP,U}(\mbf{x})$ with $x_{i} \geq \frac{N}{p^{t}}$. Then, by construction $y_{i}  = x_{i} - \frac{N}{p^{t}}$. Hence, by Observation \ref{obs_one}
\begin{align*}
C_{U}(\mbf{y}) \leq C_{U}(\mbf{x}) - \frac{N}{p^{t}} < \card{U} - 1.
\end{align*}
Thus, $\mbf{y} \in P_{\infty}$.

\item Let $0 < \epsilon < 1$, $\mbf{x}$ a vertex of $P_{t}$ where $\sum_{e \in E(G)} (X_{e} + Y_{e})$ attains a maximum, $t$ and $\mbf{y} \in P_{\infty}$ to be as in part \eqref{item_aux_lem_1} of this Lemma. Since $\norm{\mbf{y}}_{1} \leq \LP_{G,\Q}$ and $\norm{\mbf{x} - \mbf{y}}_{1} < \epsilon$, it follows that
\begin{align*}
\LP_{G,\Q,t} = \norm{\mbf{x}}_{1} \leq \norm{\mbf{x} - \mbf{y}}_{1} + \norm{\mbf{y}}_{1} \leq \epsilon + \LP_{G,\Q}.
\end{align*}
Consequently, 
\begin{align*}
\LP_{G,\Q,t} - \LP_{G,\Q} < \epsilon.
\end{align*}
This proves that $\lim_{t\ra \infty} \LP_{G,\Q,t} = \LP_{G,\Q}$.
\end{enumerate}
\eproof

\bcor
\label{cor_lct_leq_LP_G,Q}
Let $G$ be a graph, then
\begin{align*}
\lct(J_{G}) \leq \LP_{G,\Q}.
\end{align*}
\ecor

\bproof
Clear from Proposition \ref{prop_ft_leq_LP_G_Q} and Theorem \ref{thm_ft_leq_lct_and_lct_equals_lim_ft}.
\eproof

\brem
One can always bound the F-threshold of an ideal in a polynomial ring with respect to the maximal ideal by the optimal value of a linear programming problem by removing the constraints of the form $\tilde{C}_{U} \leq \card{U}-1$ from \eqref{alg_primal_LP_G} (cf. \cite{shibuta2009logcanonical} or \cite{hernandez2014Fpure}).
\erem

The previous results of this section can be summarized as
\begin{equation}
\label{eqn_chain_of_relations}
\text{binomial grade } J_{G} = r(G) = \LP_{G,\Z} \leq \ft(J_{G}) \leq \LP_{G,\Q} = \LP_{G,\Q}^{\ast} \leq \LP_{G,\Z}^{\ast} = \hgt(J_{G}).
\end{equation}

\brem
The relation $\ft(I) \leq \hgt(I)$ holds more generally, cf. \cite[Proposition 2.6(1)]{takagi2004fpure}.
\erem

\section{Proof of Conjecture \ref{conjectue_fpt_equal_LP} for Block Graphs}

\label{sec_proof_cjec_block_graphs}

\bdefn
A graph $G$ is called {\it biconnected} if $G$ is connected and $G \setminus v$ is connected for every vertex $v$ of $G$. A graph $G$ is called a {\it block graph} if every maximal biconnected component of $G$ is a complete graph. We call a maximal biconnected component of $G$ a {\it block}. A vertex $v$ of $G$ is called a {\it cut vertex} if the number of connected components of $G\setminus v$ is strictly larger than the number of connected components of $G$. 

Let $G$ be a block graph and $B$ a block of $G$, then we define $\cut(B)$ to be the set of all cut vertices of $G$ contained in $B$. We say that $B$ is a {\it leaf} of $G$ if $\card{\cut(B)} = 1$, and we say that $G$ has the {\it two block property} if every cut vertex of $G$ is adjacent to exactly two blocks of $G$. An edge of $G$ which is also a leaf of $G$ is called a \textit{whisker} of $G$.
\edefn

Throughout this section we will let $n_{G} := \card{V(G)}$, $\ell_{G}$ will denote the number of leaves of $G$, $w_{G}$ will denote the number of connected components of $G$ isomorphic to a complete graph. When the graph $G$ is clear from context we will omit the subscript $G$. 

\bprop
\label{prop_upper_bound_LP_block_graph}
Let $G$ be a block graph. Then, 
\begin{align*}
\LP_{G,\Q} \leq n_{G} - \frac{1}{2} \ell - w.
\end{align*}
\eprop

\bproof
By Remark \ref{rem_LP_G,Q_reduces_connected_cases} we may assume that $G$ is connected. If $G$ is a complete graph there is nothing to prove. We reduce to the case that $G$ is a connected block graph which is not a complete graph.

Since $G$ is a connected block graph which is not a complete graph $w = 0$ and $\ell \geq 1$. Partition the blocks of $G$ into $\mathcal{B}_{1}$, the leaves of $G$, and $\mathcal{B}_{2}$, the remaining blocks. Write $\mathcal{B}_{1} = \{B_{1},\ldots,B_{\ell}\}$. Each $B_{i} \in \mathcal{B}_{1}$ has exactly one cut vertex which we will call $c_{i}$ and at least one other vertex, i.e. $\card{V(B_{i})} \geq 2$. Let
\begin{align*}
T := \left( \bigcup_{B \in \mathcal{B}_{2}} V(B) \right) \cup \left( \bigcup \{c_{i} \mid 1 \leq i \leq \ell\} \right).
\end{align*}
For $v \in V$, define
\begin{equation}
b_{v,x} = b_{v,y} =
\begin{cases}
1/2, & v \in T \\
0, & v \in V(G) \setminus T.
\end{cases}
\end{equation}
For $U \in \A$, define
\begin{equation}
c_{U} = 
\begin{cases}
1/2, & U \in \mathcal{B}_{1} \\
1/2, & U = B_{i} \setminus \{c_{i}\} \text{ for some } 1 \leq i \leq \ell, \\
0, & \text{ else}.
\end{cases}
\end{equation}
Then, the tuple $( \{b_{v,x}\}_{v\in V}, \{b_{v,y}\}_{v\in V}, \{c_{U}\}_{U \in \A})$ is a feasible solution to the linear program \eqref{alg_dual_LP_G}, and the cost of this feasible solution is 
\begin{align*}
n_{G} - \frac{1}{2}\ell.
\end{align*}
Thus, 
\begin{align*}
\LP_{G,\Q}^{\ast} \leq n_{G} - \frac{1}{2} \ell.
\end{align*}
Hence, the result follows by strong duality, Theorem \ref{thm_strong_duality}.
\eproof

\subsection{Triangle Paths and the Two Block Property}

\bdefn
\label{defn_triangle_graph}
We say that a graph $G$ is a {\it triangle path} if every connected component of $G$ which is not an isolated vertex satisfies
\begin{enumerate}
\item the only cycles in $G$ are $K_{3}$'s which we label as $C_{1}, \ldots, C_{t}$ (we will sometimes refer to the $C_{i}$'s as the {\it triangles of} $\mathit{G}$),
\item For all $v \in G$,
\begin{align*}
\deg_{G}(v) = 
\begin{cases}
3, & v \in V(C_{i}) \text{ for some } 1 \leq i \leq t \\
1 \text{ or } 2, & \text{otherwise}.
\end{cases}
\end{align*}
\item Every vertex of $G$ belongs to at most one triangle of $G$.
\end{enumerate}

When $G$ is a triangle path and $G$ is a subgraph of $L$, then we call $G$ a {\it triangle path packing} of $L$.
\edefn

\bex

\label{ex_triangle_path}
The graph in Figure \ref{fig_ex_triangle_path} is an instance of a triangle path.
\def\putBelow{-.5}
\def\edgeDist{1.75}

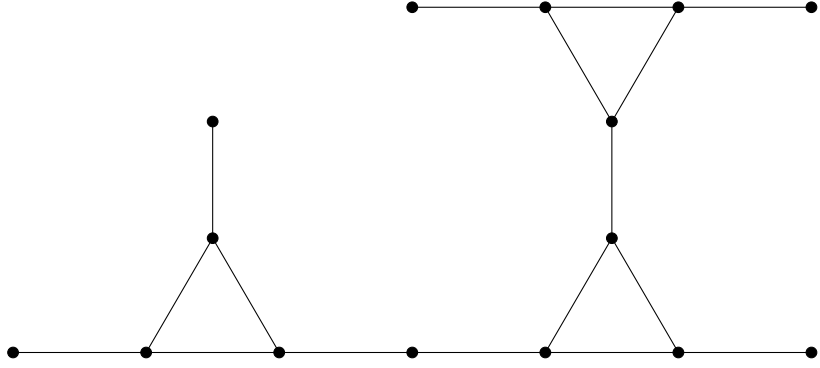
\begin{figure}[h!]
\begin{center}
\begin{tikzpicture} 
\filldraw[black] (-3*\edgeDist,0) circle (2pt) node at (-3*\edgeDist,\putBelow) {};
\filldraw[black] (-2*\edgeDist,0) circle (2pt) node at (-2*\edgeDist,\putBelow) {};
\filldraw[black] (-1*\edgeDist,0) circle (2pt) node at (-1*\edgeDist,\putBelow) {};
\filldraw[black] (0,0) circle (2pt) node at (0,\putBelow) {};
\filldraw[black] (1*\edgeDist,0) circle (2pt) node at (1*\edgeDist,\putBelow) {};
\filldraw[black] (2*\edgeDist,0) circle (2pt) node at (2*\edgeDist,\putBelow) {};
\filldraw[black] (3*\edgeDist,0) circle (2pt) node at (3*\edgeDist,\putBelow) {};
\filldraw[black] (-3/2*\edgeDist, 
0.866025403784*\edgeDist) circle (2pt) node at (-3/2*\edgeDist+0.5, 
0.866025403784*\edgeDist) {};
\filldraw[black] (3/2*\edgeDist, 
0.866025403784*\edgeDist) circle (2pt) node at (3/2*\edgeDist+0.5, 
0.866025403784*\edgeDist) {};
\filldraw[black] (-3/2*\edgeDist,1.75*\edgeDist) circle (2pt) node at (-3/2*\edgeDist+0.5,1.75*\edgeDist) {};
\filldraw[black] (3/2*\edgeDist,1.75*\edgeDist) circle (2pt) node at (3/2*\edgeDist+0.5,1.75*\edgeDist) {};
\filldraw[black] (0,2.61602540379*\edgeDist) circle (2pt) node at (0,2.61602540379*\edgeDist+.5) {};
\filldraw[black] (1*\edgeDist,2.61602540379*\edgeDist) circle (2pt) node at (1*\edgeDist,2.61602540379*\edgeDist+.5) {};
\filldraw[black] (2*\edgeDist,2.61602540379*\edgeDist) circle (2pt) node at (2*\edgeDist,2.61602540379*\edgeDist+.5) {};
\filldraw[black] (3*\edgeDist,2.61602540379*\edgeDist) circle (2pt) node at (3*\edgeDist,2.61602540379*\edgeDist+.5) {};

\draw (-3*\edgeDist,0) to (-2*\edgeDist,0);
\draw (-2*\edgeDist,0) to (-1*\edgeDist,0);
\draw (-1*\edgeDist,0) to (0,0);
\draw (-0,0) to (1*\edgeDist,0);
\draw (1*\edgeDist,0) to (2*\edgeDist,0);
\draw (2*\edgeDist,0) to (3*\edgeDist,0);
\draw (-2*\edgeDist,0) to (-3/2*\edgeDist, 
0.866025403784*\edgeDist);
\draw (-1*\edgeDist,0) to (-3/2*\edgeDist, 
0.866025403784*\edgeDist);
\draw (2*\edgeDist,0) to (3/2*\edgeDist, 
0.866025403784*\edgeDist);
\draw (1*\edgeDist,0) to (3/2*\edgeDist, 
0.866025403784*\edgeDist);
\draw (-3/2*\edgeDist,1.75*\edgeDist) to (-3/2*\edgeDist, 
0.866025403784*\edgeDist);
\draw (3/2*\edgeDist,1.75*\edgeDist) to (3/2*\edgeDist, 
0.866025403784*\edgeDist);

\draw (1*\edgeDist, 
2.61602540379*\edgeDist) to (3/2*\edgeDist,1.75*\edgeDist);
\draw (3/2*\edgeDist,1.75*\edgeDist) to (2*\edgeDist, 
2.61602540379*\edgeDist);
\draw (-0,2.61602540379*\edgeDist) to (1*\edgeDist,2.61602540379*\edgeDist);
\draw (1*\edgeDist,2.61602540379*\edgeDist) to (2*\edgeDist,2.61602540379*\edgeDist);
\draw (2*\edgeDist,2.61602540379*\edgeDist) to (3*\edgeDist,2.61602540379*\edgeDist);

\end{tikzpicture}
\end{center}
\caption{Triangle Path}
\label{fig_ex_triangle_path}
\end{figure}
\eex

\brem
\label{rem_triangle_path_has_two_block_prop}
A triangle path is an example of a family of block graphs with the two block property. 
\erem

\bprop
\label{prop_formula_for_fpt_of_triangle_path}
Let $G$ be a triangle path. Label the $K_{3}$'s in the triangle path $C_{1},\ldots, C_{t}$. Label the connected components of $G \setminus \bigcup_{k=1}^{t} E(C_{k})$ by $P_{1},\ldots, P_{s}$. Then,
\begin{equation}
\label{eqn_ft_formula_triangle_1}
\ft(J_{G}) = \sum_{i=1}^{s} \card{E(P_{i})} + \frac{3}{2}t.
\end{equation}
Moreover, 
\begin{equation}
\label{eqn_ft_formula_triangle_2}
\ft(J_{G}) = n_{G} - \frac{1}{2}\ell_{G} - w_{G}.
\end{equation}
\eprop

\bproof
By induction on $t$ it can be shown that \eqref{eqn_ft_formula_triangle_1} and \eqref{eqn_ft_formula_triangle_2} are equivalent. By Lemma \ref{lem_ft_additive} and Remark \ref{rem_ft_of_complete_graph} we may assume that $G$ is connected and not a complete graph, i.e. that $w = 0$. By Proposition \ref{prop_ft_leq_LP_G_Q} and Proposition \ref{prop_upper_bound_LP_block_graph}, and the above mentioned equivalence we have that $\ft(J_{G}) \leq \sum_{i=1}^{s} \card{E(P_{i})} + \frac{3}{2}t$. It thus suffices to prove the reverse inequality.

Let $R$ and $\m$ be as in Definition \ref{defn_bin_edge_ideal}. First consider the case that $\fchar(\K) = p > 2$. The first step in the proof is to label the vertices of the graph and to assign to every edge $e$ a tuple $(a_{e},b_{e})$ so that the following conditions are satisfied.
\begin{enumerate}
\item Some vertex having degree one has been assigned the label $1$.
\item If $a \in V(P_{j})$ and $b,c \in V(P_{i})$ for some $1 \leq i, j \leq \ell$ with the distance from $b$ to $a$ smaller than the distance from $c$ to $a$, then either $a < b  < c$ or $c < b < a$. \label{item_trianle_path_labeling_2}
\item $(a_{e},b_{e})$ equals $(0,1)$, $(1,0)$, or $(0.5,0.5)$ if $e$ is not an edge of any triangle of $G$.
\item $(a_{e},b_{e})$ equals $(0,0.5)$ or $(0.5,0)$ if $e$ is the edge of any triangle of $G$.
\item $\sum_{e \in E(G)} \left( a_{e} + b_{e} \right) = \sum_{i=1}^{s} \card{E(P_{i})} + \frac{3}{2}t$. \label{item_trianle_path_labeling_5}
\item The monomial 
\begin{equation}
\label{eqn_monomial_construction}
m :=
\prod_{\substack{e = \{i,j\} \in E(G)\\ i<j}} (X_{i} Y_{j})^{a_{ij} \cdot (p^{e}-1)} \cdot (X_{j} Y_{i})^{b_{ij} \cdot (p^{e}-1)}
\end{equation}
does not belong to the ideal $\m^{[p^{e}]}$ and appears in the polynomial
\begin{equation}
\label{eqn_polynomial_construction}
f := \prod_{\substack{e = \{i,j\} \in E(G)\\ i<j}} f_{ij}^{(a_{ij} + b_{ij})\cdot (p^{e} -1)}.
\end{equation}
with a coefficient which is non-zero mod $p$. \label{item_trianle_path_labeling_6}
\end{enumerate}
It would follow from items \eqref{item_trianle_path_labeling_5} and \eqref{item_trianle_path_labeling_6} that 
\begin{align*}
\bigg(\sum_{i=1}^{s} \card{E(P_{i})} + \frac{3}{2}t \bigg) (p^{e}-1) \leq \nu_{e}(J_{G}),
\end{align*}
and consequently that
\begin{align*}
\sum_{i=1}^{s} \card{E(P_{i})} + \frac{3}{2}t \leq \ft(J_{G}).
\end{align*}

We prove the existence of this labeling by induction on $t$. It $t = 0$, then $G$ is a path. Starting at a vertex of degree one assign it the label $1$, then label the vertices consecutively. Assign to every edge the tuple $(1,0)$. That the above conditions are satisfied is clear.

Suppose now that $t \geq 1$. There exists a triangle of $G$, call it $T$, so that removing the edges of $T$ from $G$ yields a graph having three connected components---two of which are paths, say $P$ and $Q$, and the other is a triangle path, call it $G^{'}$, having one less triangle than $G$. Let $b$ be the vertex belonging to $V(T) \cap V(G^{'})$. Let $a$ be the vertex adjacent to $b$ in $G^{'}$. Apply the induction hypothesis to $G^{'}$ and we may assume that $b$ has not been given the label $1$. By item \ref{item_trianle_path_labeling_2}, it follows that $a < b$. Label the remaining vertices of $T$ by $n+1$ and $n+2$. Label the vertex of $V(P) \cup V(Q)$ adjacent to $n+1$ (respectively $n+2$) by $n+3$ (respectively $n+4$). Label the remaining vertices of $V(P) \cup V(Q)$ so as to satisfy condition \eqref{item_trianle_path_labeling_2}. Use Table \ref{table_char_greater_two} to assign tuples to the edges of $T$ and to the edges of $P$ and $Q$ incident to the vertices of $T$ based upon the tuple assigned to the edge $\{a,b\}$. Now, $P$ and $Q$ each have an edge which has been assigned a tuple; extend this label to the remaining edges of $P$ and $Q$. After performing this construction the coefficient on $m$ in $f$ from item \eqref{item_trianle_path_labeling_6} is $\binom{p^{e}-1}{\frac{p^{e}-1}{2}}^{g}$ where $g$ is the number of times that the label $(0.5,0.5)$ appears on an edge. This binomial is non-zero by Lucas's theorem, Theorem \ref{thm_Lucas}.

\begin{table}
\begin{center}
\begin{tabular}{c|c|c|c|c|c}
$\{a,b\}$ & $\{b,n+1\}$ & $\{b,n+2\}$ & $\{n+1,n+2\}$ & $\{n+1,n+3\}$ & $\{n+2,n+4\}$ \\
\hline
$(0,1)$ & $(0,0.5)$ & $(0,0.5)$ & $(0.5,0)$ & $(0,1)$ & $(0.5,0.5)$ \\ \hline
$(1,0)$ & $(0.5,0)$ & $(0.5,0)$ & $(0,0.5)$ & $(0,1)$ & $(0.5,0.5)$ \\ \hline
$(0.5,0.5)$ & $(0.5,0)$ & $(0,0.5)$ & $(0,0.5)$ & $(1,0)$ & $(0,1)$
\end{tabular}
\end{center}
\caption{\label{table_char_greater_two}$\fchar(\K) = p > 2$}
\end{table}

Suppose now that $\fchar(\K) = 2$ and fix $r \geq 0$. We modify the above construction by using Table \ref{table_char_equal_two} in the construction. In this case we have that 
\begin{align*}
\bigg(2 \sum_{i=1}^{s} \card{E(P_{i})} + 3 t \bigg) (p^{r-1}-1) \leq \nu_{r}(J_{G}).
\end{align*}
Dividing by $p^{r}$ and taking limit as $r \ra \infty$ proves the desired inequality when $\fchar(\K) = p = 2$.	

\begin{table}
\begin{center}
\begin{tabular}{c|c|c|c|c|c}
$\{a,b\}$ & $\{b,n+1\}$ & $\{b,n+2\}$ & $\{n+1,n+2\}$ & $\{n+1,n+3\}$ & $\{n+2,n+4\}$ \\
\hline
$(0,p^{r}-1)$ & $(0,p^{r-1}-1)$ & $(0,p^{r-1}-1)$ & $(p^{r-1}-1,0)$ & $(0,p^{r}-1)$ & $(p^{r-1},p^{r-1}-1)$ \\ \hline
$(p^{r}-1,0)$ & $(p^{r-1}-1,0)$ & $(p^{r-1}-1,0)$ & $(0,p^{r-1}-1)$ & $(0,p^{r}-1)$ & $(p^{r-1},p^{r-1}-1)$ \\ \hline
$(p^{r-1},p^{r-1}-1)$ & $(p^{r-1}-1,0)$ & $(0,p^{r-1}-1)$ & $(0,p^{r-1}-1)$ & $(p^{r}-1,0)$ & $(0,p^{r}-1)$
\end{tabular}
\end{center}
\caption{\label{table_char_equal_two}$\fchar(\K) = p = 2$}
\end{table}
\eproof

\bcor
\label{fpt_splits_along_triangle_subgraphs}
Let $G$ be a triangle path and $G_{1},\ldots, G_{r}$ be subgraphs of $G$ satisfying
\begin{enumerate}
\item the $G_{i}$ are triangle paths,
\item $E(G) = \bigsqcup\limits_{i=1}^{r} E(G_{i})$,
\item every triangle of $G$ is a subgraph of some $G_{i}$
\end{enumerate}
Then,
\begin{align*}
\ft(J_{G}) = \sum_{i=1}^{r} \ft(J_{G_{i}}).
\end{align*}
\ecor

\bproof
Follows from Proposition \ref{prop_formula_for_fpt_of_triangle_path}.
\eproof

\bcor
\label{cor_ft_deleting_vertices}
Let $H$ be a triangle path. Suppose $C \subset V(H)$ satisfies that $H \setminus C$ is a triangle path, no two vertices of $C$ are adjacent via an edge of $H$, and every vertex of $C$ has degree two in $H$. Then we have
\begin{align*}
\ft(J_{H \setminus C}) = \ft(J_{H}) - 2 \card{C}.
\end{align*}
\ecor

\bproof
Since $H$ and $H \setminus C$ are triangle paths, we can compute their F-threshold as $3/2$ the number of triangles plus the number of edges not in a triangle. The later two conditions on $C$ imply that the number of triangles of $H$ and $H\setminus C$ are the same and the number of edges in $H$ not belonging to any triangle path is equal to the number of edges of $H \setminus C$ not belonging to any triangle path plus $2 \card{C}$.
\eproof

We record the following lemma which we will need later which can be viewed as an analogue of Corollary \ref{cor_ft_deleting_vertices}.

\blem
\label{lem_LP_deleting_vertices}
Let $G$ be a graph, $C \subset V(G)$, and $\{G_{i}\}_{i=1}^{r}$ the connected components of $G \setminus C$, then
\begin{align*}
\LP_{G,\Q} \leq \sum_{i=1}^{r} \LP_{G_{i},\Q} + 2 \card{C}.
\end{align*}
\elem

\bproof
By strong duality, Theorem \ref{thm_strong_duality}, it suffices to prove the statement for the dual linear program. For $i =1,\ldots,r$, let $(\{b_{i,v,x}\}_{v\in V}, \{b_{i,v,y}\}_{v\in V}, \{c_{i,U}\}_{U\in\A_{i}})$ be an optimal solution to the linear program \eqref{alg_dual_LP_G} realizing $\LP^{\ast}_{G_{i},\Q}$. These {\it local} solutions together with the assignments $b_{v,x} := 1$, $b_{v,y} := 1$ for $v \in C$, and $c_{U} := 0$ for $U \in \A \setminus \cup_{i=1}^{r} \A_{i}$ yield a feasible solution to the linear program \eqref{alg_dual_LP_G} realizing $\LP^{\ast}_{G,\Q}$.
\eproof

The next proposition shows that the F-threshold of a block graph with the two block property can be computed by a \textit{maximal} triangle path packing.

\bprop
\label{prop_good_triangle_path_packing}
Let $G$ be a block graph with the two block property. Let $\ell$ denote the number of leaves of $G$ and $w$ the number of connected components of $G$ which are isomorphic to a complete graph. 
Then, $G$ contains a subgraph $H$ satisfying 
\begin{enumerate}
\item $H$ is a triangle path packing, \label{prop_1}
\item $V(G) = V(H)$, \label{prop_2}
\item for every block $B$ of $G$, $B$ contains at most one triangle of $H$, \label{prop_3}
\item a block $B$ of $G$ contains a triangle of $H$ only if $\card{\cut{B}}$ is odd, \label{prop_3.5}
\item \label{prop_4} if $T$ is a triangle of $H$ contained in a block $B$ of $G$, then $T$ contains two cut vertices unless $\cut(B) = V(B)$ in which case $T$ contains three, 
\item a block $B$ of $G$ contains exactly one vertex of $H$ of degree one if and only if $B$ is a leaf of $G$, \label{prop_5}
\item a block $B$ of $G$ contains a vertex of $H$ of degree zero if and only if $B$ is a single vertex, \label{prop_5.5}
\item if $T$ is a triangle of $H$ contained in $B$, $c \in \cut(B) \cap V(T)$, and $\tilde{B} \neq B$ is the block of $G$ incident to $c$, then the edge of $H$ incident to $c$ and not contained in $T$ is contained in $\tilde{B}$, \label{prop_6}
\item $\ft(J_{H}) = n - \frac{1}{2}\ell - w$.
\label{prop_7}
\end{enumerate}
\eprop

The proof idea is to induce on the number of blocks of $G$. We start by choosing a block of $G$ which is not a leaf and so that every incident block with the exception of at most one is a leaf of $G$. In Figure \ref{fig_block_graph}, blocks $\{12,13\}$ and $\{1,2,3,4,5,6\}$ are the only such blocks. We will consider block $B := \{1,2,3,4,5,6\}$. Next, we choose two leaves of $B$ to delete (if $B$ has only one leaf then delete that leaf), Figure \ref{fig_graph_deleting_leaves}. By induction we construct a triangle path packing satisfying the above conditions, Figure \ref{fig_triangle_path_subgraph}. Lastly, extend the triangle path packing of the subgraph to a triangle path packing of the original graph so that the new triangle path packing has the desired properties, Figure \ref{fig_triangle_path_original_graph}. When the number of cut vertices of $G$ in $B$ is less than or equal to $3$ several cases must be considered to extend the triangle path packing while preserving the above properties. However, when the number of cut vertices of $G$ in $B$ is larger than $3$ then the triangle path packing can be extended by taking a Hamiltonian path on the deleted vertices.

\def\edgeDistEx{.9}
\def\edgeDistExHor{.75}
\def\putAbove{.35}

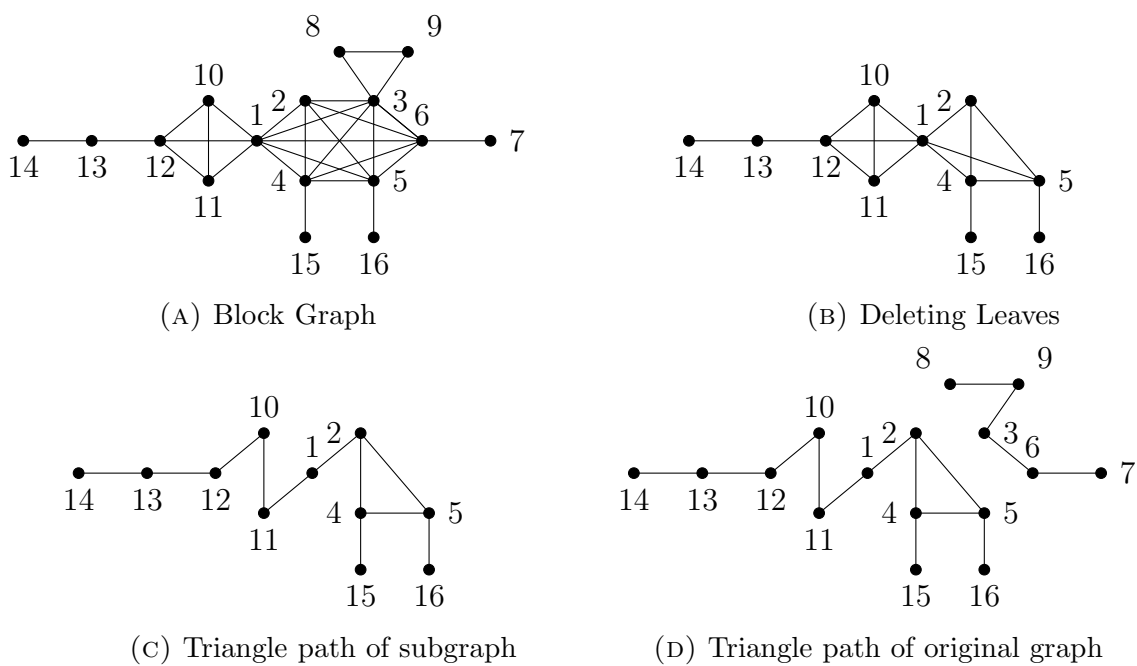
\begin{figure}[ht]
\begin{subfigure}{.45\linewidth}

\begin{tikzpicture} 

\filldraw[black] (0*\edgeDistEx,0) circle (2pt) node at (0*\edgeDistEx,\putAbove) {$1$};

\filldraw[black] (.707*\edgeDistEx,0.707*\edgeDistExHor) circle (2pt) node at (.707*\edgeDistEx-\putAbove,0.707*\edgeDistExHor) {$2$};

\filldraw[black] (1.707*\edgeDistEx,.707*\edgeDistExHor) circle (2pt) node at (1.707*\edgeDistEx+\putAbove,.707*\edgeDistExHor) {$3$};

\filldraw[black] (.707*\edgeDistEx,-0.707*\edgeDistExHor) circle (2pt) node at (.707*\edgeDistEx-\putAbove,-0.707*\edgeDistExHor) {$4$};

\filldraw[black] (1.707*\edgeDistEx,-0.707*\edgeDistExHor) circle (2pt) node at (1.707*\edgeDistEx+\putAbove,-0.707*\edgeDistExHor) {$5$};

\filldraw[black] (2.414*\edgeDistEx,0*\edgeDistExHor) circle (2pt) node at (2.414*\edgeDistEx,0*\edgeDistExHor+\putAbove) {$6$};

\filldraw[black] (3.414*\edgeDistEx,0*\edgeDistExHor) circle (2pt) node at (3.414*\edgeDistEx+\putAbove,0*\edgeDistExHor) {$7$};

\filldraw[black] (1.207*\edgeDistEx,1.573*\edgeDistExHor) circle (2pt) node at (1.207*\edgeDistEx-\putAbove,1.573*\edgeDistExHor+\putAbove) {$8$};

\filldraw[black] (2.207*\edgeDistEx,1.573*\edgeDistExHor) circle (2pt) node at (2.207*\edgeDistEx+\putAbove,1.573*\edgeDistExHor+\putAbove) {$9$};

\filldraw[black] (-1.414*\edgeDistEx,0*\edgeDistExHor) circle (2pt) node at (-1.414*\edgeDistEx,0*\edgeDistExHor-\putAbove) {$12$};

\filldraw[black] (-2.414*\edgeDistEx,0*\edgeDistExHor) circle (2pt) node at (-2.414*\edgeDistEx,0*\edgeDistExHor-\putAbove) {$13$};

\filldraw[black] (-3.414*\edgeDistEx,0*\edgeDistExHor) circle (2pt) node at (-3.414*\edgeDistEx,0*\edgeDistExHor-\putAbove) {$14$};

\filldraw[black] (-.707*\edgeDistEx,.707*\edgeDistExHor) circle (2pt) node at (-.707*\edgeDistEx,.707*\edgeDistExHor+\putAbove) {$10$};

\filldraw[black] (-.707*\edgeDistEx,-.707*\edgeDistExHor) circle (2pt) node at (-.707*\edgeDistEx,-.707*\edgeDistExHor-\putAbove) {$11$};

\filldraw[black] (.707*\edgeDistEx,-1.707*\edgeDistExHor) circle (2pt) node at (.707*\edgeDistEx,-1.707*\edgeDistExHor-\putAbove) {$15$};

\filldraw[black] (1.707*\edgeDistEx,-1.707*\edgeDistExHor) circle (2pt) node at (1.707*\edgeDistEx,-1.707*\edgeDistExHor-\putAbove) {$16$};

\draw (-3.414*\edgeDistEx,0*\edgeDistExHor) to (-2.414*\edgeDistEx,0*\edgeDistExHor); 

\draw (-2.414*\edgeDistEx,0*\edgeDistExHor) to  (-1.414*\edgeDistEx,0*\edgeDistExHor);
\draw (-1.414*\edgeDistEx,0*\edgeDistExHor) to  (0*\edgeDistEx,0); 
\draw (-1.414*\edgeDistEx,0*\edgeDistExHor) to  (-.707*\edgeDistEx,.707*\edgeDistExHor); 
\draw (-1.414*\edgeDistEx,0*\edgeDistExHor) to  (-.707*\edgeDistEx,-.707*\edgeDistExHor); 
\draw (0*\edgeDistEx,0) to  (.707*\edgeDistEx,0.707*\edgeDistExHor); 
\draw (0*\edgeDistEx,0) to  (1.707*\edgeDistEx,.707*\edgeDistExHor); 
\draw  (0*\edgeDistEx,0) to  (.707*\edgeDistEx,-0.707*\edgeDistExHor);
\draw (0*\edgeDistEx,0) to  (1.707*\edgeDistEx,-0.707*\edgeDistExHor);
\draw (0*\edgeDistEx,0) to (2.414*\edgeDistEx,0*\edgeDistExHor); 
\draw  (.707*\edgeDistEx,0.707*\edgeDistExHor) to  (1.707*\edgeDistEx,.707*\edgeDistExHor); 
\draw (.707*\edgeDistEx,0.707*\edgeDistExHor) to  (.707*\edgeDistEx,-0.707*\edgeDistExHor);
\draw (.707*\edgeDistEx,0.707*\edgeDistExHor) to (1.707*\edgeDistEx,-0.707*\edgeDistExHor);
\draw (.707*\edgeDistEx,0.707*\edgeDistExHor) to (2.414*\edgeDistEx,0*\edgeDistExHor) ;
\draw  (1.707*\edgeDistEx,.707*\edgeDistExHor) to  (.707*\edgeDistEx,-0.707*\edgeDistExHor);
\draw  (1.707*\edgeDistEx,.707*\edgeDistExHor) to  (1.707*\edgeDistEx,-0.707*\edgeDistExHor);
\draw (1.707*\edgeDistEx,.707*\edgeDistExHor) to (2.414*\edgeDistEx,0*\edgeDistExHor) ;
\draw (.707*\edgeDistEx,-0.707*\edgeDistExHor) to  (1.707*\edgeDistEx,-0.707*\edgeDistExHor);
\draw (.707*\edgeDistEx,-0.707*\edgeDistExHor) to  (2.414*\edgeDistEx,0*\edgeDistExHor);
\draw (1.707*\edgeDistEx,-0.707*\edgeDistExHor) to  (2.414*\edgeDistEx,0*\edgeDistExHor);
\draw (2.414*\edgeDistEx,0*\edgeDistExHor) to (3.414*\edgeDistEx,0*\edgeDistExHor); 
\draw (1.707*\edgeDistEx,-0.707*\edgeDistExHor) to  (1.707*\edgeDistEx,-1.707*\edgeDistExHor);
\draw (1.707*\edgeDistEx,.707*\edgeDistExHor) to  (1.207*\edgeDistEx,1.573*\edgeDistExHor);
\draw (1.707*\edgeDistEx,.707*\edgeDistExHor) to (2.207*\edgeDistEx,1.573*\edgeDistExHor);
\draw (1.207*\edgeDistEx,1.573*\edgeDistExHor) to (2.207*\edgeDistEx,1.573*\edgeDistExHor);
\draw (1.707*\edgeDistEx,.707*\edgeDistExHor) to (2.414*\edgeDistEx,0*\edgeDistExHor) ;
\draw (.707*\edgeDistEx,-0.707*\edgeDistExHor) to (.707*\edgeDistEx,-1.707*\edgeDistExHor);
\draw (-.707*\edgeDistEx,.707*\edgeDistExHor) to  (0*\edgeDistEx,0) ; 
\draw (-.707*\edgeDistEx,.707*\edgeDistExHor) to  (-.707*\edgeDistEx,-.707*\edgeDistExHor) ; 
\draw (0*\edgeDistEx,0) to (-.707*\edgeDistEx,-.707*\edgeDistExHor); 

\end{tikzpicture}
\subcaption{Block Graph}
\label{fig_block_graph}
\end{subfigure}
\hfill
%
\begin{subfigure}{.45\linewidth}
\begin{tikzpicture}

\filldraw[black] (0*\edgeDistEx,0) circle (2pt) node at (0*\edgeDistEx,\putAbove) {$1$};

\filldraw[black] (.707*\edgeDistEx,0.707*\edgeDistExHor) circle (2pt) node at (.707*\edgeDistEx-\putAbove,0.707*\edgeDistExHor) {$2$};


\filldraw[black] (.707*\edgeDistEx,-0.707*\edgeDistExHor) circle (2pt) node at (.707*\edgeDistEx-\putAbove,-0.707*\edgeDistExHor) {$4$};

\filldraw[black] (1.707*\edgeDistEx,-0.707*\edgeDistExHor) circle (2pt) node at (1.707*\edgeDistEx+\putAbove,-0.707*\edgeDistExHor) {$5$};





\filldraw[black] (-1.414*\edgeDistEx,0*\edgeDistExHor) circle (2pt) node at (-1.414*\edgeDistEx,0*\edgeDistExHor-\putAbove) {$12$};

\filldraw[black] (-2.414*\edgeDistEx,0*\edgeDistExHor) circle (2pt) node at (-2.414*\edgeDistEx,0*\edgeDistExHor-\putAbove) {$13$};

\filldraw[black] (-3.414*\edgeDistEx,0*\edgeDistExHor) circle (2pt) node at (-3.414*\edgeDistEx,0*\edgeDistExHor-\putAbove) {$14$};

\filldraw[black] (-.707*\edgeDistEx,.707*\edgeDistExHor) circle (2pt) node at (-.707*\edgeDistEx,.707*\edgeDistExHor+\putAbove) {$10$};

\filldraw[black] (-.707*\edgeDistEx,-.707*\edgeDistExHor) circle (2pt) node at (-.707*\edgeDistEx,-.707*\edgeDistExHor-\putAbove) {$11$};

\filldraw[black] (.707*\edgeDistEx,-1.707*\edgeDistExHor) circle (2pt) node at (.707*\edgeDistEx,-1.707*\edgeDistExHor-\putAbove) {$15$};

\filldraw[black] (1.707*\edgeDistEx,-1.707*\edgeDistExHor) circle (2pt) node at (1.707*\edgeDistEx,-1.707*\edgeDistExHor-\putAbove) {$16$};

\draw (-3.414*\edgeDistEx,0*\edgeDistExHor) to (-2.414*\edgeDistEx,0*\edgeDistExHor); 

\draw (-2.414*\edgeDistEx,0*\edgeDistExHor) to  (-1.414*\edgeDistEx,0*\edgeDistExHor);
\draw (-1.414*\edgeDistEx,0*\edgeDistExHor) to  (0*\edgeDistEx,0); 
\draw (-1.414*\edgeDistEx,0*\edgeDistExHor) to  (-.707*\edgeDistEx,.707*\edgeDistExHor); 
\draw (-1.414*\edgeDistEx,0*\edgeDistExHor) to  (-.707*\edgeDistEx,-.707*\edgeDistExHor); 
\draw (0*\edgeDistEx,0) to  (.707*\edgeDistEx,0.707*\edgeDistExHor); 
\draw  (0*\edgeDistEx,0) to  (.707*\edgeDistEx,-0.707*\edgeDistExHor);
\draw (0*\edgeDistEx,0) to  (1.707*\edgeDistEx,-0.707*\edgeDistExHor);
\draw (.707*\edgeDistEx,0.707*\edgeDistExHor) to  (.707*\edgeDistEx,-0.707*\edgeDistExHor);
\draw (.707*\edgeDistEx,0.707*\edgeDistExHor) to (1.707*\edgeDistEx,-0.707*\edgeDistExHor);
\draw (.707*\edgeDistEx,-0.707*\edgeDistExHor) to  (1.707*\edgeDistEx,-0.707*\edgeDistExHor);
\draw (1.707*\edgeDistEx,-0.707*\edgeDistExHor) to  (1.707*\edgeDistEx,-1.707*\edgeDistExHor);
\draw (.707*\edgeDistEx,-0.707*\edgeDistExHor) to (.707*\edgeDistEx,-1.707*\edgeDistExHor);
\draw (-.707*\edgeDistEx,.707*\edgeDistExHor) to  (0*\edgeDistEx,0) ; 
\draw (-.707*\edgeDistEx,.707*\edgeDistExHor) to  (-.707*\edgeDistEx,-.707*\edgeDistExHor) ; 
\draw (0*\edgeDistEx,0) to (-.707*\edgeDistEx,-.707*\edgeDistExHor); 

\end{tikzpicture}

\subcaption{Deleting Leaves}
\label{fig_graph_deleting_leaves}
\end{subfigure}
\hfill
%
\begin{subfigure}{.45\linewidth}
\begin{tikzpicture}

\filldraw[black] (0*\edgeDistEx,0) circle (2pt) node at (0*\edgeDistEx,\putAbove) {$1$};

\filldraw[black] (.707*\edgeDistEx,0.707*\edgeDistExHor) circle (2pt) node at (.707*\edgeDistEx-\putAbove,0.707*\edgeDistExHor) {$2$};


\filldraw[black] (.707*\edgeDistEx,-0.707*\edgeDistExHor) circle (2pt) node at (.707*\edgeDistEx-\putAbove,-0.707*\edgeDistExHor) {$4$};

\filldraw[black] (1.707*\edgeDistEx,-0.707*\edgeDistExHor) circle (2pt) node at (1.707*\edgeDistEx+\putAbove,-0.707*\edgeDistExHor) {$5$};





\filldraw[black] (-1.414*\edgeDistEx,0*\edgeDistExHor) circle (2pt) node at (-1.414*\edgeDistEx,0*\edgeDistExHor-\putAbove) {$12$};

\filldraw[black] (-2.414*\edgeDistEx,0*\edgeDistExHor) circle (2pt) node at (-2.414*\edgeDistEx,0*\edgeDistExHor-\putAbove) {$13$};

\filldraw[black] (-3.414*\edgeDistEx,0*\edgeDistExHor) circle (2pt) node at (-3.414*\edgeDistEx,0*\edgeDistExHor-\putAbove) {$14$};

\filldraw[black] (-.707*\edgeDistEx,.707*\edgeDistExHor) circle (2pt) node at (-.707*\edgeDistEx,.707*\edgeDistExHor+\putAbove) {$10$};

\filldraw[black] (-.707*\edgeDistEx,-.707*\edgeDistExHor) circle (2pt) node at (-.707*\edgeDistEx,-.707*\edgeDistExHor-\putAbove) {$11$};

\filldraw[black] (.707*\edgeDistEx,-1.707*\edgeDistExHor) circle (2pt) node at (.707*\edgeDistEx,-1.707*\edgeDistExHor-\putAbove) {$15$};

\filldraw[black] (1.707*\edgeDistEx,-1.707*\edgeDistExHor) circle (2pt) node at (1.707*\edgeDistEx,-1.707*\edgeDistExHor-\putAbove) {$16$};

\draw (-3.414*\edgeDistEx,0*\edgeDistExHor) to (-2.414*\edgeDistEx,0*\edgeDistExHor); 

\draw (-2.414*\edgeDistEx,0*\edgeDistExHor) to  (-1.414*\edgeDistEx,0*\edgeDistExHor);
\draw (-1.414*\edgeDistEx,0*\edgeDistExHor) to  (-.707*\edgeDistEx,.707*\edgeDistExHor); 
\draw (0*\edgeDistEx,0) to  (.707*\edgeDistEx,0.707*\edgeDistExHor); 
\draw (.707*\edgeDistEx,0.707*\edgeDistExHor) to  (.707*\edgeDistEx,-0.707*\edgeDistExHor);
\draw (.707*\edgeDistEx,0.707*\edgeDistExHor) to (1.707*\edgeDistEx,-0.707*\edgeDistExHor);
\draw (.707*\edgeDistEx,-0.707*\edgeDistExHor) to  (1.707*\edgeDistEx,-0.707*\edgeDistExHor);
\draw (1.707*\edgeDistEx,-0.707*\edgeDistExHor) to  (1.707*\edgeDistEx,-1.707*\edgeDistExHor);
\draw (.707*\edgeDistEx,-0.707*\edgeDistExHor) to (.707*\edgeDistEx,-1.707*\edgeDistExHor);
\draw (-.707*\edgeDistEx,.707*\edgeDistExHor) to  (-.707*\edgeDistEx,-.707*\edgeDistExHor) ; 
\draw (0*\edgeDistEx,0) to (-.707*\edgeDistEx,-.707*\edgeDistExHor); 

\end{tikzpicture}
\subcaption{Triangle path of subgraph}
\label{fig_triangle_path_subgraph}
\end{subfigure}
\begin{subfigure}{.45\linewidth}
\begin{tikzpicture}

\filldraw[black] (0*\edgeDistEx,0) circle (2pt) node at (0*\edgeDistEx,\putAbove) {$1$};

\filldraw[black] (.707*\edgeDistEx,0.707*\edgeDistExHor) circle (2pt) node at (.707*\edgeDistEx-\putAbove,0.707*\edgeDistExHor) {$2$};

\filldraw[black] (1.707*\edgeDistEx,.707*\edgeDistExHor) circle (2pt) node at (1.707*\edgeDistEx+\putAbove,.707*\edgeDistExHor) {$3$};

\filldraw[black] (.707*\edgeDistEx,-0.707*\edgeDistExHor) circle (2pt) node at (.707*\edgeDistEx-\putAbove,-0.707*\edgeDistExHor) {$4$};

\filldraw[black] (1.707*\edgeDistEx,-0.707*\edgeDistExHor) circle (2pt) node at (1.707*\edgeDistEx+\putAbove,-0.707*\edgeDistExHor) {$5$};

\filldraw[black] (2.414*\edgeDistEx,0*\edgeDistExHor) circle (2pt) node at (2.414*\edgeDistEx,0*\edgeDistExHor+\putAbove) {$6$};

\filldraw[black] (3.414*\edgeDistEx,0*\edgeDistExHor) circle (2pt) node at (3.414*\edgeDistEx+\putAbove,0*\edgeDistExHor) {$7$};

\filldraw[black] (1.207*\edgeDistEx,1.573*\edgeDistExHor) circle (2pt) node at (1.207*\edgeDistEx-\putAbove,1.573*\edgeDistExHor+\putAbove) {$8$};

\filldraw[black] (2.207*\edgeDistEx,1.573*\edgeDistExHor) circle (2pt) node at (2.207*\edgeDistEx+\putAbove,1.573*\edgeDistExHor+\putAbove) {$9$};

\filldraw[black] (-1.414*\edgeDistEx,0*\edgeDistExHor) circle (2pt) node at (-1.414*\edgeDistEx,0*\edgeDistExHor-\putAbove) {$12$};

\filldraw[black] (-2.414*\edgeDistEx,0*\edgeDistExHor) circle (2pt) node at (-2.414*\edgeDistEx,0*\edgeDistExHor-\putAbove) {$13$};

\filldraw[black] (-3.414*\edgeDistEx,0*\edgeDistExHor) circle (2pt) node at (-3.414*\edgeDistEx,0*\edgeDistExHor-\putAbove) {$14$};

\filldraw[black] (-.707*\edgeDistEx,.707*\edgeDistExHor) circle (2pt) node at (-.707*\edgeDistEx,.707*\edgeDistExHor+\putAbove) {$10$};

\filldraw[black] (-.707*\edgeDistEx,-.707*\edgeDistExHor) circle (2pt) node at (-.707*\edgeDistEx,-.707*\edgeDistExHor-\putAbove) {$11$};

\filldraw[black] (.707*\edgeDistEx,-1.707*\edgeDistExHor) circle (2pt) node at (.707*\edgeDistEx,-1.707*\edgeDistExHor-\putAbove) {$15$};

\filldraw[black] (1.707*\edgeDistEx,-1.707*\edgeDistExHor) circle (2pt) node at (1.707*\edgeDistEx,-1.707*\edgeDistExHor-\putAbove) {$16$};

\draw (-3.414*\edgeDistEx,0*\edgeDistExHor) to (-2.414*\edgeDistEx,0*\edgeDistExHor); 

\draw (-2.414*\edgeDistEx,0*\edgeDistExHor) to  (-1.414*\edgeDistEx,0*\edgeDistExHor);
\draw (-1.414*\edgeDistEx,0*\edgeDistExHor) to  (-.707*\edgeDistEx,.707*\edgeDistExHor); 
\draw (0*\edgeDistEx,0) to  (.707*\edgeDistEx,0.707*\edgeDistExHor); 
\draw (.707*\edgeDistEx,0.707*\edgeDistExHor) to  (.707*\edgeDistEx,-0.707*\edgeDistExHor);
\draw (.707*\edgeDistEx,0.707*\edgeDistExHor) to (1.707*\edgeDistEx,-0.707*\edgeDistExHor);
\draw (1.707*\edgeDistEx,.707*\edgeDistExHor) to (2.414*\edgeDistEx,0*\edgeDistExHor) ;
\draw (.707*\edgeDistEx,-0.707*\edgeDistExHor) to  (1.707*\edgeDistEx,-0.707*\edgeDistExHor);
\draw (2.414*\edgeDistEx,0*\edgeDistExHor) to (3.414*\edgeDistEx,0*\edgeDistExHor); 
\draw (1.707*\edgeDistEx,-0.707*\edgeDistExHor) to  (1.707*\edgeDistEx,-1.707*\edgeDistExHor);
\draw (1.707*\edgeDistEx,.707*\edgeDistExHor) to (2.207*\edgeDistEx,1.573*\edgeDistExHor);
\draw (1.207*\edgeDistEx,1.573*\edgeDistExHor) to (2.207*\edgeDistEx,1.573*\edgeDistExHor);
\draw (.707*\edgeDistEx,-0.707*\edgeDistExHor) to (.707*\edgeDistEx,-1.707*\edgeDistExHor);
\draw (-.707*\edgeDistEx,.707*\edgeDistExHor) to  (-.707*\edgeDistEx,-.707*\edgeDistExHor) ; 
\draw (0*\edgeDistEx,0) to (-.707*\edgeDistEx,-.707*\edgeDistExHor); 

\end{tikzpicture}
\subcaption{Triangle path of original graph}
\label{fig_triangle_path_original_graph}
\end{subfigure}
\caption{Proof Idea of Proposition \ref{prop_good_triangle_path_packing}}
\label{fig_picture_illustrating_proof}
\end{figure}

\bproof
It is enough to show this proposition for each connected component of $G$. If $G$ is isomorphic to a complete graph, then take $H$ to be a Hamiltonian path and the Proposition follows since $\ell = 0$ and $w = 1$ in this case. Thus, we may assume henceforth that $G$ is a connected graph and not isomorphic to a complete graph.

The proof is by induction on the number of blocks of $G$ which we denote by $t$. Since $G$ is not a complete graph, $t \geq 2$. If $t = 2$, then $G$ is two complete graphs joined at a cut vertex. In this case, $\ell = 2$ and $G$ has a Hamiltonian path. Take $H$ to be a Hamiltonian path of $G$, and it is clear that $H$ satisfies the conditions listed above.

Suppose now that $t \geq 3$ and that the Proposition has been shown for all block graphs having the two block property with at most $t-1$ blocks. Let $B$ be a block of $G$ which is not a leaf such that every block incident to $B$ with the exception of at most one is a leaf of $G$.  Let $r$ denote the number of cut vertices of $B$. Observe that $r \geq 2$ since $B$ is not a leaf.

\ul{Case $r = 2$:} Let $L$ be a leaf of $B$. Let $\tilde{c}$ be the cut vertex of $B$ which is adjacent to $L$. If $\card{V(B)} > 2$ put $\tilde{G} := G \setminus L $; otherwise, put $\tilde{G} := G \setminus (L \setminus \tilde{c})$. Let $\tilde{B}$ be the image of $B$ in $\tilde{G}$. Since $\tilde{G}$ has $t-1$ blocks, by induction there exists a triangle path $\tilde{H}$ of $\tilde{G}$ satisfying the above properties. Since $\tilde{B}$ is a leaf in $\tilde{G}$, $\tilde{H} \mid_{\tilde{B}}$ has a vertex of degree one, say $v$, by property \eqref{prop_5}. Consider now $\tilde{H}$ as a subgraph of $G$, and construct $H$ a triangle path in $G$ as follows. Observe that if $\card{V(B)} = 2$, then $v = \tilde{c}$. Otherwise, adjoin an edge from $v$ to $\tilde{c}$. In either case $H$ is now a subgraph of $G$ containing $\tilde{c}$ as a degree one vertex. Finally, adjoin to $H$ a Hamiltonian path contained in $L$ which has a terminal vertex at $\tilde{c}$. Now, $H$ is a triangle path of $G$ which satisfies properties \eqref{prop_1} through \eqref{prop_7}. Moreover, $\ft(J_{H})$ is of the desired form by Proposition \ref{prop_formula_for_fpt_of_triangle_path}, induction hypothesis, and the fact that we obtained $H$ from $\tilde{H}$ by adding a path.

\ul{Case $r = 3$:} Let $\cut{B} = \{c_{1},c_{2},c_{3}\}$. Let $L_{i}$ be leaves adjacent to $B$ at $c_{i}$ for $i \in \{1,2\}$.

\ul{Sub-case $\card{V(B)} = 3$:} Let $\tilde{G}$ be the graph attained after removing $L_{1}$, $L_{2}$, and $B \setminus \{c_{3}\}$. Note that $\tilde{G}$ is adjacent to $B$ at $c_{3}$ and that $\tilde{G}$ may or may not be a leaf. Consider the graph $G^{'}$ which is $\tilde{G}$ adjoined a whisker at the vertex $c_{3}$. Denote by $n^{'}$ and $\ell^{'}$ the number of vertices and leaves of $G^{'}$, respectively. By induction hypothesis there exists a triangle path $H^{'}$ for $G^{'}$ satisfying the above properties and 
\begin{align*}
\ft (J_{H^{'}}) = n^{'} - \frac{1}{2} \ell^{'}.
\end{align*}
First, consider the sub-sub-case that there exists a triangle in $H^{'}$ containing $c_{3}$; label the other vertices of this triangle as $a$ and $b$. Construct $H$ from $H^{'}$ by deleting the whisker attached at $c_{3}$, deleting the edge $\{a,b\}$, adding the edge $\{c_{1},c_{2}\}$, and adding Hamiltonian paths in $L_{1}$ and $L_{2}$ which have a terminal vertex at $c_{1}$ and $c_{2}$, respectively. Observe that $H$ is a triangle path, and by Proposition \ref{prop_formula_for_fpt_of_triangle_path} we compute that 
\begin{align*}
\ft(J_{H}) = \ft(J_{H^{'}}) + (-1 - \frac{1}{2} + 1) + (\card{V(L_{1})}-1) + 1 + \card{V(L_{2})}-1)
\end{align*}
where the term $(-1 - \frac{1}{2} + 1)$ is counting in $H^{'}$ the deletion of a whisker, the deletion of an edge of a triangle, and that two edges of a triangle are now counted as edges belonging to a path. The desired formula on $\ft(J_{H})$ now follows after simplifying and using that $n = n^{'} + (\card{V(L_{1})} -1) + \card{V(L_{2})}$ and that $\ell = \ell^{'} + 1$.

Second, consider the sub-sub-case that $c_{3}$ has degree two in $H^{'}$; call the neighbor of $c_{3}$ contained in $\tilde{G}$ to be $a$. Construct $H$ from $H^{'}$ by deleting the whisker, adding the edges $\{c_{1},c_{2}\}$, $\{c_{2},c_{3}\}$, $\{c_{1},c_{3}\}$ and adding Hamiltonian paths in $L_{1}$ and $L_{2}$ which have a terminal vertex at $c_{1}$ and $c_{2}$, respectively. Observe that $H$ is a triangle path, and by Proposition \ref{prop_formula_for_fpt_of_triangle_path} we compute that 
\begin{align*}
\ft(J_{H}) = \big( \ft(J_{H^{'}}) -1 \big) + \frac{3}{2} + (\card{V(L_{1})}-1)  + \card{V(L_{2})}-1)
\end{align*}
where the $-1$ is the deletion of the whisker and the plus $\frac{3}{2}$ is the addition of the triangle. The formula for $\ft(J_{H})$ now follows after simplifying and using that $n = n^{'} + (\card{V(L_{1})}-1) + \card{V(L_{2})}$ and that $\ell = \ell^{'} + 1$.

\ul{Sub-case $\card{V(B)} > 3$:} Let $\tilde{G} := G \setminus ( L_{1} \cup L_{2})$. By induction hypothesis there exists a triangle path $\tilde{H}$ of $\tilde{G}$ satisfying the desired properties. Since $\tilde{B}$ is a leaf in $\tilde{G}$, property \eqref{prop_5} guarantees the existence of a vertex $v$ in $\tilde{H} \cap \tilde{B}$ having degree one. Consider $\tilde{H}$ as a subgraph of $G$, and construct $H$ a triangle path of $G$ from $\tilde{H}$ by adding the edges $\{v,c_{1}\}$, $\{v,c_{2}\}$, $\{c_{1},c_{2}\}$ and adding Hamiltonian paths in $L_{1}$ and $L_{2}$ which have a terminal vertex at $c_{1}$ and $c_{2}$, respectively. Observe that $H$ is a triangle path, and by Proposition \ref{prop_formula_for_fpt_of_triangle_path} we compute that 
\begin{align*}
\ft(J_{H}) = \ft(J_{\tilde{H}})  + \frac{3}{2} + (\card{V(L_{1})}-1) + (\card{V(L_{2})}-1)
\end{align*}
where the plus $\frac{3}{2}$ is the addition of the triangle. The desired formula on $\ft(J_{H})$ now follows after simplifying and using that $n = \tilde{n} + \card{V(L_{1})} + \card{V(L_{2})}$ and that $\ell = \tilde{\ell} + 1$.

\ul{Case $r > 3$:} Let $c_{1}$ and $c_{2}$ be two cut vertices of $B$ which have leaves $L_{1}$ and $L_{2}$, respectively. Let $\tilde{G} := G \setminus ( L_{1} \cup L_{2})$. By induction hypothesis there exists a triangle path $\tilde{H}$ of $\tilde{G}$ satisfying the desired properties. Consider $\tilde{H}$ as a subgraph of $G$, and construct $H$ a triangle path of $G$ from $\tilde{H}$ by adding the edges $\{c_{1},c_{2}\}$ and adding Hamiltonian paths in $L_{1}$ and $L_{2}$ which have a terminal vertex at $c_{1}$ and $c_{2}$, respectively. Observe that $H$ is a triangle path, and by Proposition \ref{prop_formula_for_fpt_of_triangle_path} we compute that 
\begin{align*}
\ft(J_{H}) = \ft(J_{\tilde{H}}) + (\card{V(L_{1})}-1) + 1 + \card{V(L_{2})}-1)
\end{align*}
where the plus $1$ is the addition of the edge $\{c_{1},c_{2}\}$. The desired formula on $\ft(J_{H})$ now follows after simplifying and using that $n = \tilde{n} + \card{V(L_{1})} + \card{V(L_{2})}$ and that $\ell = \tilde{\ell} + 2$.
\eproof

\bcor
\label{cor_ft_two_block_graph_equals_ft_triangle_path_packing}
Let $G$ be a block graph having the two block property. Let $\ell$ denote the number of leaves of $G$ and $w$ the number of connected components of $G$ isomorphic to a complete graph. Then, there exists a triangle path packing $H$ of $G$ such that
\begin{align*}
\ft(J_{H}) = \ft(J_{G}) = \LP_{G,\Q} = n - \frac{1}{2} \ell - w.
\end{align*}
\ecor

\bproof
Follows from Lemma \ref{lem_ft_monotonic}, Proposition \ref{prop_ft_leq_LP_G_Q}, Proposition \ref{prop_good_triangle_path_packing}, and Proposition \ref{prop_upper_bound_LP_block_graph}.
\eproof

\subsection{Two Block Approximation and General Case}

Next, we show how to reduce the computation of the F-threshold of a block graph to that of a block graph with the two block property.

\bdefn
Let $G$ be a block graph. We call $H$ a {\it two block approximation} of $G$ if $H$ can be obtained via the following process
\begin{enumerate}
\item Enumerate the cut vertices of $G$ incident to at least three blocks of $G$; call this set $C$.
\item For each $c \in C$, choose two blocks incident to $c$, say $A_{c}$ and $B_{c}$.
\item Delete all edges incident to $c$ and not contained in $E(A_{c}) \cup E(B_{c})$.
\end{enumerate}
\edefn

\def\edgeDistBlockApproxX{.75}
\def\edgeDistBlockApproxY{.75}

\begin{figure}[ht]

\begin{subfigure}{.3\linewidth}
\begin{center}
\begin{tikzpicture} 
\filldraw (-0.707106*\edgeDistBlockApproxX,-0.707106*\edgeDistBlockApproxY) circle (2pt); 
\filldraw (0*\edgeDistBlockApproxX,0*\edgeDistBlockApproxY) circle (2pt); 
\filldraw (1*\edgeDistBlockApproxX,0*\edgeDistBlockApproxY) circle (2pt); 
\filldraw (1*\edgeDistBlockApproxX,1*\edgeDistBlockApproxY) circle (2pt); 
\filldraw (0*\edgeDistBlockApproxX,1*\edgeDistBlockApproxY) circle (2pt); 
\filldraw (1*\edgeDistBlockApproxX,-1*\edgeDistBlockApproxY) circle (2pt); 
\filldraw (1.707106*\edgeDistBlockApproxX, -0.707106*\edgeDistBlockApproxY) circle (2pt); 
\filldraw (2*\edgeDistBlockApproxX,0*\edgeDistBlockApproxY) circle (2pt); 
\filldraw (-1*\edgeDistBlockApproxX,1*\edgeDistBlockApproxY) circle (2pt); 
\filldraw (-0.707106*\edgeDistBlockApproxX,1.707106*\edgeDistBlockApproxY) circle (2pt); 
\filldraw (0*\edgeDistBlockApproxX,2*\edgeDistBlockApproxY) circle (2pt); 
\filldraw (1.707106*\edgeDistBlockApproxX,1.707106*\edgeDistBlockApproxY) circle (2pt); 

\draw (-0.707106*\edgeDistBlockApproxX,-0.707106*\edgeDistBlockApproxY) to (0,0); 
\draw (0,0) to (1*\edgeDistBlockApproxX,0*\edgeDistBlockApproxY) ; 
\draw (0,0) to  (1*\edgeDistBlockApproxX,1*\edgeDistBlockApproxY) ; 
\draw (0,0) to  (0*\edgeDistBlockApproxX,1*\edgeDistBlockApproxY) ; 

\draw (1*\edgeDistBlockApproxX,0) to  (1*\edgeDistBlockApproxX,1*\edgeDistBlockApproxY) ; 
\draw (1*\edgeDistBlockApproxX,0) to  (0*\edgeDistBlockApproxX,1*\edgeDistBlockApproxY) ; 
\draw (1*\edgeDistBlockApproxX,0) to  (1*\edgeDistBlockApproxX,-1*\edgeDistBlockApproxY) ; 
\draw (1*\edgeDistBlockApproxX,0) to  (1.707106*\edgeDistBlockApproxX, -0.707106*\edgeDistBlockApproxY) ; 
\draw (1*\edgeDistBlockApproxX,0) to  (2*\edgeDistBlockApproxX,0*\edgeDistBlockApproxY) ; 

\draw (1*\edgeDistBlockApproxX,1*\edgeDistBlockApproxY) to  (0*\edgeDistBlockApproxX,1*\edgeDistBlockApproxY) ; 
\draw (1*\edgeDistBlockApproxX,1*\edgeDistBlockApproxY) to  (1.707106*\edgeDistBlockApproxX,1.707106*\edgeDistBlockApproxY) ; 

\draw (0*\edgeDistBlockApproxX,1*\edgeDistBlockApproxY) to  (-1*\edgeDistBlockApproxX,1*\edgeDistBlockApproxY) ; 
\draw (0*\edgeDistBlockApproxX,1*\edgeDistBlockApproxY) to  (-0.707106*\edgeDistBlockApproxX,1.707106*\edgeDistBlockApproxY) ; 
\draw (0*\edgeDistBlockApproxX,1*\edgeDistBlockApproxY) to  (0*\edgeDistBlockApproxX,2*\edgeDistBlockApproxY) ; 

\draw (1.707106*\edgeDistBlockApproxX, -0.707106*\edgeDistBlockApproxY) to  (2*\edgeDistBlockApproxX,0*\edgeDistBlockApproxY) ; 

\end{tikzpicture}
\end{center}
\subcaption{Block Graph}
\label{fig_block_graph_too_many_blocks}
\end{subfigure}
\hfill
\begin{subfigure}{.3\linewidth}
\begin{center}
\begin{tikzpicture} 
\filldraw (-0.707106*\edgeDistBlockApproxX,-0.707106*\edgeDistBlockApproxY) circle (2pt); 
\filldraw (0*\edgeDistBlockApproxX,0*\edgeDistBlockApproxY) circle (2pt); 
\filldraw (1*\edgeDistBlockApproxX,0*\edgeDistBlockApproxY) circle (2pt); 
\filldraw (1*\edgeDistBlockApproxX,1*\edgeDistBlockApproxY) circle (2pt); 
\filldraw (0*\edgeDistBlockApproxX,1*\edgeDistBlockApproxY) circle (2pt); 
\filldraw (1*\edgeDistBlockApproxX,-1*\edgeDistBlockApproxY) circle (2pt); 
\filldraw (1.707106*\edgeDistBlockApproxX, -0.707106*\edgeDistBlockApproxY) circle (2pt); 
\filldraw (2*\edgeDistBlockApproxX,0*\edgeDistBlockApproxY) circle (2pt); 
\filldraw (-1*\edgeDistBlockApproxX,1*\edgeDistBlockApproxY) circle (2pt); 
\filldraw (-0.707106*\edgeDistBlockApproxX,1.707106*\edgeDistBlockApproxY) circle (2pt); 
\filldraw (0*\edgeDistBlockApproxX,2*\edgeDistBlockApproxY) circle (2pt); 
\filldraw (1.707106*\edgeDistBlockApproxX,1.707106*\edgeDistBlockApproxY) circle (2pt); 

\draw (-0.707106*\edgeDistBlockApproxX,-0.707106*\edgeDistBlockApproxY) to (0,0); 
\draw (0,0) to  (1*\edgeDistBlockApproxX,1*\edgeDistBlockApproxY) ; 

\draw (1*\edgeDistBlockApproxX,0) to  (1*\edgeDistBlockApproxX,-1*\edgeDistBlockApproxY) ; 
\draw (1*\edgeDistBlockApproxX,0) to  (1.707106*\edgeDistBlockApproxX, -0.707106*\edgeDistBlockApproxY) ; 
\draw (1*\edgeDistBlockApproxX,0) to  (2*\edgeDistBlockApproxX,0*\edgeDistBlockApproxY) ; 

\draw (1*\edgeDistBlockApproxX,1*\edgeDistBlockApproxY) to  (1.707106*\edgeDistBlockApproxX,1.707106*\edgeDistBlockApproxY) ; 

\draw (0*\edgeDistBlockApproxX,1*\edgeDistBlockApproxY) to  (-1*\edgeDistBlockApproxX,1*\edgeDistBlockApproxY) ; 
\draw (0*\edgeDistBlockApproxX,1*\edgeDistBlockApproxY) to  (-0.707106*\edgeDistBlockApproxX,1.707106*\edgeDistBlockApproxY) ; 

\draw (1.707106*\edgeDistBlockApproxX, -0.707106*\edgeDistBlockApproxY) to  (2*\edgeDistBlockApproxX,0*\edgeDistBlockApproxY) ; 

\end{tikzpicture}
\end{center}
\subcaption{Two Block Approximation}
\label{fig_two_block_approx_1}
\end{subfigure}
\hfill
%
\begin{subfigure}{.3\linewidth}
\begin{center}
\begin{tikzpicture} 
\filldraw (-0.707106*\edgeDistBlockApproxX,-0.707106*\edgeDistBlockApproxY) circle (2pt); 
\filldraw (0*\edgeDistBlockApproxX,0*\edgeDistBlockApproxY) circle (2pt); 
\filldraw (1*\edgeDistBlockApproxX,0*\edgeDistBlockApproxY) circle (2pt); 
\filldraw (1*\edgeDistBlockApproxX,1*\edgeDistBlockApproxY) circle (2pt); 
\filldraw (0*\edgeDistBlockApproxX,1*\edgeDistBlockApproxY) circle (2pt); 
\filldraw (1*\edgeDistBlockApproxX,-1*\edgeDistBlockApproxY) circle (2pt); 
\filldraw (1.707106*\edgeDistBlockApproxX, -0.707106*\edgeDistBlockApproxY) circle (2pt); 
\filldraw (2*\edgeDistBlockApproxX,0*\edgeDistBlockApproxY) circle (2pt); 
\filldraw (-1*\edgeDistBlockApproxX,1*\edgeDistBlockApproxY) circle (2pt); 
\filldraw (-0.707106*\edgeDistBlockApproxX,1.707106*\edgeDistBlockApproxY) circle (2pt); 
\filldraw (0*\edgeDistBlockApproxX,2*\edgeDistBlockApproxY) circle (2pt); 
\filldraw (1.707106*\edgeDistBlockApproxX,1.707106*\edgeDistBlockApproxY) circle (2pt); 

\draw (-0.707106*\edgeDistBlockApproxX,-0.707106*\edgeDistBlockApproxY) to (0,0); 
\draw (0,0) to (1*\edgeDistBlockApproxX,0*\edgeDistBlockApproxY) ; 
\draw (0,0) to  (1*\edgeDistBlockApproxX,1*\edgeDistBlockApproxY) ; 
\draw (0,0) to  (0*\edgeDistBlockApproxX,1*\edgeDistBlockApproxY) ; 

\draw (1*\edgeDistBlockApproxX,0) to  (1*\edgeDistBlockApproxX,1*\edgeDistBlockApproxY) ; 
\draw (1*\edgeDistBlockApproxX,0) to  (0*\edgeDistBlockApproxX,1*\edgeDistBlockApproxY) ; 
\draw (1*\edgeDistBlockApproxX,0) to  (1*\edgeDistBlockApproxX,-1*\edgeDistBlockApproxY) ; 

\draw (1*\edgeDistBlockApproxX,1*\edgeDistBlockApproxY) to  (0*\edgeDistBlockApproxX,1*\edgeDistBlockApproxY) ; 
\draw (1*\edgeDistBlockApproxX,1*\edgeDistBlockApproxY) to  (1.707106*\edgeDistBlockApproxX,1.707106*\edgeDistBlockApproxY) ; 

\draw (0*\edgeDistBlockApproxX,1*\edgeDistBlockApproxY) to  (-0.707106*\edgeDistBlockApproxX,1.707106*\edgeDistBlockApproxY) ; 

\draw (1.707106*\edgeDistBlockApproxX, -0.707106*\edgeDistBlockApproxY) to  (2*\edgeDistBlockApproxX,0*\edgeDistBlockApproxY) ; 

\end{tikzpicture}
\end{center}
\subcaption{Two Block Approximation}
\label{fig_two_block_approx_2}
\end{subfigure}
\caption{Examples of Two Block Approximation}
\end{figure}
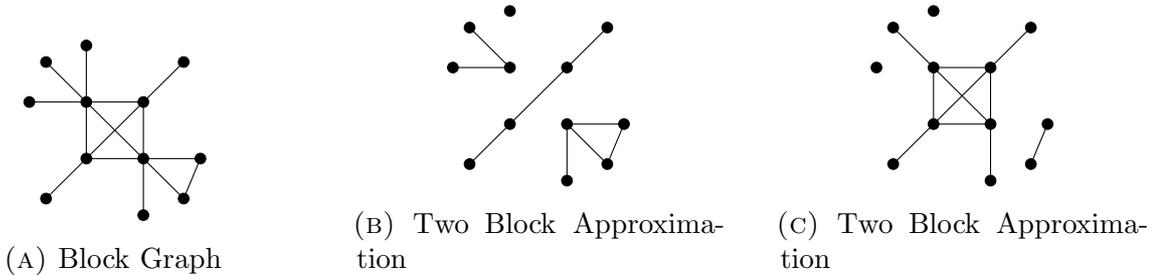

\bex
Let $G$ be the block graph in Figure \ref{fig_block_graph_too_many_blocks}. Then, the graphs in Figure \ref{fig_two_block_approx_1} and \ref{fig_two_block_approx_2} are two possible block approximations of $G$. Observe that the F-threshold of the graphs in Figure \ref{fig_two_block_approx_1} and \ref{fig_two_block_approx_2} are $8$ and $7$, respectively.

\eex

Our next goal is to show

\bthm
\label{thm_formula_ft_block_graph}
Let $G$ be a block graph, then the following quantities are equal
\begin{enumerate}
\item $\max\{ \ft J_{H} \mid \; H \text{ is a triangle path packing of } G\} \label{item_thm_1}$
\item $\max\{ \ft J_{K} \; \mid K \text{ is a two block approximation of } G\} \label{item_thm_2}$
\item $\ft(J_{G}) \label{item_thm_3}$
\item $\lct(J_{G}) \label{item_thm_4}$
\item $\max\{ \LP_{K,\Q} \mid K \text{ is a two block approximation of } G\}  \label{item_thm_5}$ \label{item_LP_K_Q}
\item $\LP_{G,\Q}  \label{item_thm_6}$.
\end{enumerate}
\ethm

\blem
\label{lem_max_ft_triangle_packing_equals_max_ft_block_approximation}
In Theorem \ref{thm_formula_ft_block_graph} items \eqref{item_thm_1}, \eqref{item_thm_2}, and \eqref{item_thm_5} are equal.
\elem

\bproof
By Corollary \ref{cor_ft_two_block_graph_equals_ft_triangle_path_packing} it follows that \eqref{item_thm_1} is greater than or equal to \eqref{item_thm_2} which is equal to \eqref{item_thm_5}. Any triangle path of $G$ has the two block property by Remark \ref{rem_triangle_path_has_two_block_prop}, and hence can be realized as a subgraph of a two block approximation of $G$. This proves that \eqref{item_thm_1} is less than or equal to \eqref{item_thm_2}. 
\eproof

To prove Theorem \ref{thm_formula_ft_block_graph} it suffices to prove that \eqref{item_thm_1} is equal to \eqref{item_thm_6}. Indeed, \eqref{item_thm_1}, \eqref{item_thm_2}, and \eqref{item_thm_5} are equal by Lemma \ref{lem_max_ft_triangle_packing_equals_max_ft_block_approximation}, and equality of \eqref{item_thm_1} and \eqref{item_thm_3} would follow from Lemma \ref{lem_ft_monotonic} and Proposition \ref{prop_ft_leq_LP_G_Q}. In particular, this shows that the computation of $\ft(J_{G})$ is independent of characteristic, and equality of \eqref{item_thm_3} and \eqref{item_thm_4} would follow from Theorem \ref{thm_ft_leq_lct_and_lct_equals_lim_ft} and Corollary \ref{cor_lct_leq_LP_G,Q}.

For the proof of Theorem \ref{thm_formula_ft_block_graph} we take the following setup.

Let $G$ be a connected block graph, $K$ a two block approximation of $G$ realizing \eqref{item_LP_K_Q} in Theorem \ref{thm_formula_ft_block_graph}, and $H$ a path packing of $K$ given by Lemma \ref{prop_good_triangle_path_packing}. Take $C$ to be the set of cut vertices of $G$ incident to at least three blocks of $G$ such that in $K$ at least one of these blocks becomes a leaf or an isolated connected component. For $c \in C$, let $B_{c,1},B_{c,2},B_{c,3},\ldots,B_{c,m_{c}}$ denote the blocks of $G$ incident to $c$ for some $m_{c} \geq 3$. Without loss of generality, we may suppose that $N_{H}(c) \subset V(B_{c,1}) \cup V(B_{c,2})$ and that $B_{c,3}$ becomes a leaf or an isolated connected component of $K$. Since $B_{c,3}$ is a leaf or isolated component of $K$, Proposition \ref{prop_good_triangle_path_packing} implies that there exists a vertex $v_{c} \in V(B_{c,3}) \cap V(H)$ having degree one or zero in $H$. 

\blem
\label{lem_C_is_empty}
With the above setup suppose that $C = \varnothing$. Then $G$ satisfies the conclusions of Theorem \ref{thm_formula_ft_block_graph}. Moreover, $\LP_{K,\Q} = \LP_{G,\Q}$.
\elem

\bproof
Because $C = \varnothing$, $\ell_{K} = \ell_{G}$ and $w_{K} = w_{G}$. Now, the result follows by Lemma \ref{lem_ft_monotonic}, Proposition \ref{prop_upper_bound_LP_block_graph}, and Corollary \ref{cor_ft_two_block_graph_equals_ft_triangle_path_packing}.
\eproof

\blem
\label{lem_vertices_incident_leaves_blocks}
With the setup from above let $c \in C$. Then
\begin{enumerate}
\item $c$ does not belong to a triangle of $H$, \label{item_vertex_not_in_triangle}
\item each vertex of $N_{H}(c)$ is not contained in a triangle of $H$, \label{item_nbrs_not_in_triangle}
\item $\card{N_{H}(c)} = 2$, and \label{item_vertex_has_two_nbrs}
\item $N_{H}(c) \not\subseteq B_{c,j}$ for $j = 1,2$, \label{item_nbrs_not_contained_single_blk}
\item no two vertices of $C$ are adjacent via an edge of $H$. \label{item_vertices_disconnected}
\end{enumerate}
\elem

\bproof
\begin{enumerate}
\item Suppose by contradiction that $c$ belongs to a triangle of $H$, i.e. that without loss of generality there exists vertices $a,b \in B_{c,1}$ so that $\{a,b\}$, $\{a,c\}$, and $\{b,c\}$ belong to $E(H)$. Let $\tilde{H}$ be the triangle path defined from $H$ by deleting the edges $\{c,a\}$ and $\{c,b\}$ and adding the edge $\{c,v_{c}\}$. Notice that no new cycles are created in $\tilde{H}$ after these modifications or else it would contradict $G$ being a block graph. Moreover,
\begin{align*}
\ft J_{\tilde{H}} = \ft J_{H} + \frac{1}{2}
\end{align*}
because in obtaining $\tilde{H}$ from $H$ we delete two edges of $H$ having weight $\frac{1}{2}$, increase the weight of one edge of $H$ from weight $\frac{1}{2}$ to weight $1$ in $\tilde{H}$, and add one edge to $\tilde{H}$ having weight $1$. This contradicts our choice of $K$ via Lemma \ref{lem_max_ft_triangle_packing_equals_max_ft_block_approximation}.

\item Suppose by contradiction that there exists $a \in N_{H}(c) \cap V(B_{c,1})$ and vertices $b$ and $d$ of $B_{c,1}$ with $\{a,b\}$, $\{a,d\}$, and $\{b,d\}$ edges of $H$. Construct $\tilde{H}$ a triangle path from $H$ by deleting the edges $\{c,a\}$ and $\{b,d\}$ and adding the edge $\{c,v_{c}\}$. Then, $\ft J_{\tilde{H}} = \frac{1}{2} + \ft J_{H}$; a contradiction.

\item By Item \eqref{item_nbrs_not_in_triangle} $\card{N_{H}(c)}$ is $1$ or $2$. Suppose by contradiction that $N_{H}(c) = \{a\} \subseteq B_{i,1}$. Construct $\tilde{H}$ from $H$ by adding the edge $\{c,v_{c}\}$. Then, $\ft(J_{\tilde{H}}) > \ft(J_{H})$; a contradiction.

\item Suppose by contradiction that $\{a,b\} = N_{H}(c) \subseteq B_{c,1}$. By the previous item, $a$ and $b$ do not belong to a triangle of $H$. If $a$ and $b$ have degree two in $H$, then define $\tilde{H}$ from $H$ by adding the edge $\{a,b\}$ and the edge $\{c,v_{c}\}$. In the case that $a$ has degree one and $b$ has degree two in $H$, define $\tilde{H}$ from $H$ by deleting the edge $\{c,b\}$ and adding the edges $\{a,b\}$ and $\{c,v_{c}\}$. The case that $b$ has degree one and $a$ has degree two in $B_{c,1}$ is handled analogously. If $a$ and $b$ both have degree one in $H$, define $\tilde{H}$ from $H$ by deleting the edge $\{a,c\}$ and adding the edges $\{a,b\}$ and $\{c,v_{c}\}$. Notice that in all cases that $\tilde{H}$ is a triangle path, a subgraph of $G$, and that $\ft(J_{\tilde{H}}) > \ft(J_{H})$; a contradiction.

\item Suppose by contradiction that $c, c^{'} \in C$ and that $\{c,c^{'}\} \in E(H)$. Construct $\tilde{H}$ a triangle path from $H$ by deleting the edge $\{c,c^{'}\}$ and adding the edges $\{c,v_{c}\}$ and $\{c^{'},v_{c^{'}}\}$. Then, $\ft(J_{\tilde{H}}) > \ft(J_{H})$; a contradiction.

\end{enumerate}
\eproof

\begin{proof}[Proof of Theorem \ref{thm_formula_ft_block_graph}]
Let the setup be as above. We induce on the number of blocks of $G$. If the number of blocks of $G$ is one, then there is nothing to prove. Suppose that the number of blocks of $G$ is strictly larger than one. If $C = \varnothing$, then we are done by Lemma \ref{lem_C_is_empty}. Suppose that $C \neq \varnothing$. Let $\{G_{i}\}_{i = 1}^{r}$ denote the connected components of $G \setminus C$. For $1 \leq i \leq r$, let $H_{i}$ denote the restriction of $H$ to $G_{i}$. We show
\begin{claim}
\label{claim_ft_H_i_equals_LP_G_i}
For $1 \leq i \leq r$, 
\begin{align*}
\label{eqn_ft_H_i_equals_LP_G_i}
\ft(J_{H_{i}}) = \LP_{G_{i},\Q}.
\end{align*}
\end{claim}

\begin{pfclaim}
By induction on number of blocks and the fact that the number of blocks of $G_{i}$ is less than the number of blocks of $G$, it suffices to show that 
\begin{align*}
    \ft(J_{H_{i}}) = \max\{ \ft J_{\tilde{H}_{i}} \mid \; \tilde{H}_{i} \text{ is a triangle path of } G_{i}\}
\end{align*}
for every $1 \leq i \leq r$. Suppose by contradiction that $\ft(J_{H_{i}}) < \ft(J_{\tilde{H}_{i}})$ for some fixed choice of $1 \leq i \leq r$ and a triangle path $\tilde{H}_{i}$ of $G_{i}$. Let $c_{1},\ldots,c_{t}$ be the vertices of $C$ which are adjacent to some vertex of $G_{i}$ via an edge of $H$; call this vertex $c_{j}^{'}$ for $1 \leq j \leq t$. Construct $\tilde{H}$ from $H$ as follows. 
\begin{enumerate}[(i)]
\item Delete the edges $\{c_{j},c_{j}^{'}\}$ for $1 \leq j \leq t$ \label{item_constr_1}.
\item Add the edges $\{c_{j},v_{c_{j}}\}$ for $1 \leq j \leq t$ \label{item_constr_2}.
\item Replace $H_{i}$ on $G_{i}$ by $\tilde{H}_{i}$ \label{item_constr_3}.
\end{enumerate}
Now, $\tilde{H}$ is a triangle path. Indeed parts \eqref{item_constr_1} and \eqref{item_constr_2} of the construction do not create any cycles in $\tilde{H}$; otherwise, $G$ could not be a block graph. After performing parts \eqref{item_constr_1} and \eqref{item_constr_2}, $H_{i}$ forms its own connected component by Lemma \ref{lem_vertices_incident_leaves_blocks} item \eqref{item_nbrs_not_contained_single_blk}. Thus after performing operation \eqref{item_constr_3} $\tilde{H}$ will remain a triangle path. Clearly, $\ft(J_{\tilde{H}}) > \ft(J_{H})$ which contradicts our choice of $H$, Lemma \ref{lem_max_ft_triangle_packing_equals_max_ft_block_approximation}.
\end{pfclaim}

Observe that
\begin{align*}
\ft(J_{H}) \leq \LP_{G,\Q} &\leq \sum_{i=1}^{r} \LP_{G_{i},\Q} +2 \card{C}  & & (\text{Lemma } \ref{lem_LP_deleting_vertices})\\
&= \sum_{i=1}^{r} \ft(J_{H_{i}}) + 2 \card{C} &  & (\text{Claim } \ref{claim_ft_H_i_equals_LP_G_i}) \\
&= \ft(J_{H}) & & (\text{Lemma } \ref{lem_vertices_incident_leaves_blocks} \text{ parts } \eqref{item_vertex_not_in_triangle}, \, \eqref{item_nbrs_not_in_triangle}, \, \eqref{item_vertex_has_two_nbrs}, \, \eqref{item_vertices_disconnected} \\
& & & \text{ and Corollary } \ref{cor_ft_deleting_vertices} ).
\end{align*}
This proves equality of \eqref{item_thm_1} and \eqref{item_thm_6} and completes the proof.
\end{proof}

\bcor
\label{cor_main_conj_true_block_graphs}
Conjecture \ref{conjectue_fpt_equal_LP} is true when $G$ is a block graph.
\ecor

\bcor
If $G$ is a tree, then
\begin{align*}
r(G) = \ft(J_{G}).
\end{align*}
\ecor

\bproof
By Theorem \ref{thm_formula_ft_block_graph}, $\ft(J_{G}) = \max\{ \ft J_{K} \; \mid K \text{ is a two block approximation of } G\}$. Any two block approximation $K$ of $G$ is a semi-path of $G$ since $G$ is a tree, and hence $\ft(J_{G}) \leq r(G)$. The inequality $r(G) \leq \ft(J_{G})$ is Equation \eqref{eqn_chain_of_relations}.
\eproof

In the next section we show the stronger statement that for a tree $G$ 
\begin{align*}
r(G) = \hgt(J_{G}).
\end{align*}

\subsection{Remark on the Case of Chordal Graphs}

The next proposition shows that if Conjecture \ref{conjectue_fpt_equal_LP} is true for chordal graphs, then the proof techniques used for block graphs will not suffice. We show:

\bprop
\label{prop_bad_chordal_graph}
There exists a chordal graph $G$ such that 
\begin{align*}
\max\{\ft(J_{H}) \mid H \subset G \text{ is triangle path packing} \} < \LP_{G,\Q}.
\end{align*}
\eprop

\bproof
Let $G$ be the graph below. Let $G_{1}$ and $G_{2}$ be the subgraphs of $G$ induced on the vertices $\{1,2,3,4,5,6,13,14\}$ and $\{7,8,9,10,11,12,15,16\}$, respectively. Let $H$ be a triangle path packing of $G$ maximal wrt $\ft(J_{H})$. Let $H_{1}$ and $H_{2}$ be the restriction of $H$ to $G_{1}$ and $G_{2}$ respectively. For $i = 1,2$, let $t_{i}$, $e_{i}$ denote the number of triangles of $H_{i}$ and edges of $H_{i}$ not belonging to any triangle, respectively.

\def\putBelow{-.5}
\def\edgeDistChor{1.25}

\begin{figure}[h!]
\begin{center}
\begin{tikzpicture} 
\filldraw (-11/2*\edgeDistChor,0) circle (2pt) node at (-11/2*\edgeDistChor,\putBelow) {$1$};
\filldraw (-9/2*\edgeDistChor,0) circle (2pt) node at (-9/2*\edgeDistChor,\putBelow) {$2$};
\filldraw (-7/2*\edgeDistChor,0) circle (2pt) node at (-7/2*\edgeDistChor,\putBelow) {$3$};
\filldraw (-5/2*\edgeDistChor,0) circle (2pt) node at (-5/2*\edgeDistChor,\putBelow) {$4$};
\filldraw (-3/2*\edgeDistChor,0) circle (2pt) node at (-3/2*\edgeDistChor,\putBelow) {$5$};
\filldraw (-1/2*\edgeDistChor,0) circle (2pt) node at (-1/2*\edgeDistChor,\putBelow) {$6$};
\filldraw (-3*\edgeDistChor, 
0.866025403784*\edgeDistChor) circle (2pt) node at (-3*\edgeDistChor+0.5, 
0.866025403784*\edgeDistChor) {$13$};
\filldraw (-3*\edgeDistChor,1.75*\edgeDistChor) circle (2pt) node at (-3*\edgeDistChor+0.5,1.75*\edgeDistChor) {$14$};

\draw (-11/2*\edgeDistChor,0) to (-9/2*\edgeDistChor,0); 
\draw (-9/2*\edgeDistChor,0) to (-7/2*\edgeDistChor,0); 
\draw (-7/2*\edgeDistChor,0) to (-5/2*\edgeDistChor,0); 
\draw (-5/2*\edgeDistChor,0) to (-3/2*\edgeDistChor,0); 
\draw (-3/2*\edgeDistChor,0) to (-1/2*\edgeDistChor,0); 
\draw (-7/2*\edgeDistChor,0) to (-3*\edgeDistChor, 
0.866025403784*\edgeDistChor); 
\draw (-5/2*\edgeDistChor,0) to (-3*\edgeDistChor, 
0.866025403784*\edgeDistChor); 
\draw (-3*\edgeDistChor,1.75*\edgeDistChor) to (-3*\edgeDistChor, 
0.866025403784*\edgeDistChor); 

\draw (-11/2*\edgeDistChor,0) .. controls (-5*\edgeDistChor,-1.5*\edgeDistChor) and (-1*\edgeDistChor,-1.5*\edgeDistChor) .. (-1/2*\edgeDistChor,0); 

\draw (-11/2*\edgeDistChor,0) .. controls (-5*\edgeDistChor,-1/3*\edgeDistChor) and (-4*\edgeDistChor,-1/3*\edgeDistChor) .. (-7/2*\edgeDistChor,0); 

\draw (-5/2*\edgeDistChor,0) .. controls (-2*\edgeDistChor,-1/3*\edgeDistChor) and (-1*\edgeDistChor,-1/3*\edgeDistChor) .. (-1/2*\edgeDistChor,0); 

\draw (-7/2*\edgeDistChor,0) .. controls (-3*\edgeDistChor,-1*\edgeDistChor) and (-1*\edgeDistChor,-1*\edgeDistChor) .. (-1/2*\edgeDistChor,0); 


\filldraw (1/2*\edgeDistChor,0) circle (2pt) node at (1/2*\edgeDistChor,\putBelow) {$7$};
\filldraw (3/2*\edgeDistChor,0) circle (2pt) node at (3/2*\edgeDistChor,\putBelow) {$8$};
\filldraw (5/2*\edgeDistChor,0) circle (2pt) node at (5/2*\edgeDistChor,\putBelow) {$9$};
\filldraw (7/2*\edgeDistChor,0) circle (2pt) node at (7/2*\edgeDistChor,\putBelow) {$10$};
\filldraw (9/2*\edgeDistChor,0) circle (2pt) node at (9/2*\edgeDistChor,\putBelow) {$11$};
\filldraw (11/2*\edgeDistChor,0) circle (2pt) node at (11/2*\edgeDistChor,\putBelow) {$12$};
\filldraw (3*\edgeDistChor, 
0.866025403784*\edgeDistChor) circle (2pt) node at (3*\edgeDistChor+0.5, 
0.866025403784*\edgeDistChor) {$15$};
\filldraw (3*\edgeDistChor,1.75*\edgeDistChor) circle (2pt) node at (3*\edgeDistChor+0.5,1.75*\edgeDistChor) {$16$};

\draw (-1/2*\edgeDistChor,0) to (1/2*\edgeDistChor,0); 

\draw (1/2*\edgeDistChor,0) to (3/2*\edgeDistChor,0); 
\draw (3/2*\edgeDistChor,0) to (5/2*\edgeDistChor,0); 
\draw (5/2*\edgeDistChor,0) to (7/2*\edgeDistChor,0); 
\draw (7/2*\edgeDistChor,0) to (9/2*\edgeDistChor,0); 
\draw (9/2*\edgeDistChor,0) to (11/2*\edgeDistChor,0); 
\draw (5/2*\edgeDistChor,0) to (3*\edgeDistChor, 
0.866025403784*\edgeDistChor); 
\draw (7/2*\edgeDistChor,0) to (3*\edgeDistChor, 
0.866025403784*\edgeDistChor); 
\draw (3*\edgeDistChor,1.75*\edgeDistChor) to (3*\edgeDistChor, 
0.866025403784*\edgeDistChor); 

\draw (1/2*\edgeDistChor,0) .. controls (1*\edgeDistChor,-1.5*\edgeDistChor) and (5*\edgeDistChor,-1.5*\edgeDistChor) .. (11/2*\edgeDistChor,0); 

\draw (1/2*\edgeDistChor,0) .. controls (1*\edgeDistChor,-1/3*\edgeDistChor) and (2*\edgeDistChor,-1/3*\edgeDistChor) .. (5/2*\edgeDistChor,0); 

\draw (7/2*\edgeDistChor,0) .. controls (4*\edgeDistChor,-1/3*\edgeDistChor) and (5*\edgeDistChor,-1/3*\edgeDistChor) .. (11/2*\edgeDistChor,0); 

\draw (1/2*\edgeDistChor,0) .. controls (1*\edgeDistChor,-1*\edgeDistChor) and (3*\edgeDistChor,-1*\edgeDistChor) .. (7/2*\edgeDistChor,0); 

\end{tikzpicture}
\end{center}
\caption{Chordal Graph}
\label{fig_bad_chordal_graph}
\end{figure}
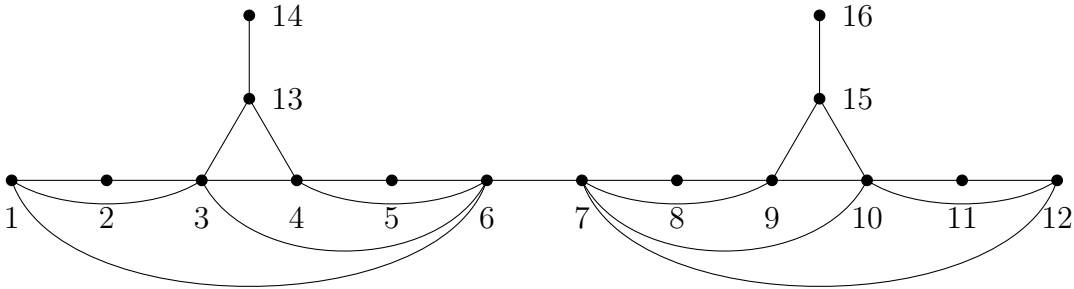

\begin{claim}
\label{clm_eq_comps}
$e_{1} + \frac{3}{2} t_{1} = e_{2} + \frac{3}{2} t_{2}$. 
\end{claim}

\begin{pfclaim}
Without loss of generality suppose $e_{1} + \frac{3}{2} t_{1} < e_{2} + \frac{3}{2} t_{2}$. Construct $\tilde{H}$ from $H$ by deleting $H_{1}$ and replicating $H_{2}$ on $V(G_{1})$. Then, $\tilde{H}$ is a triangle path packing with $\ft(J_{H}) < \ft(J_{\tilde{H}})$; a contradiction. 
\end{pfclaim}

It follows from Proposition \ref{prop_formula_for_fpt_of_triangle_path} and Claim \ref{clm_eq_comps} that $\ft(J_{H})$ is an integer. However, $\LP_{G,\Q} = 14.5$ which can be realized by the following assignment
\begin{align*}
(a_{[1, 2]}, b_{[1, 2]}) &= (0.0,1.0) & (a_{[9, 10]}, b_{[9, 10]}) &= (0.0,0.5) & (a_{[4, 13]}, b_{[4, 13]}) &= (0.0,0.5)\\
(a_{[2, 3]}, b_{[2, 3]}) &= (0.0,1.0) & (a_{[10, 11]}, b_{[10, 11]}) &= (0.0,0.5) & (a_{[7, 9]}, b_{[7, 9]}) &= (0.0,0.0)\\
(a_{[3, 4]}, b_{[3, 4]}) &= (0.0,0.5) & (a_{[11, 12]}, b_{[11, 12]}) &= (0.5,0.5) & (a_{[7, 10]}, b_{[7, 10]}) &= (0.0,0.0)\\
(a_{[4, 5]}, b_{[4, 5]}) &= (0.5,0.5) & (a_{[1, 3]}, b_{[1, 3]}) &= (0.0,0.0) & (a_{[7, 12]}, b_{[7, 12]}) &= (0.5,0.0)\\
(a_{[5, 6]}, b_{[5, 6]}) &= (0.0,0.5) & (a_{[1, 6]}, b_{[1, 6]}) &= (1.0,0.0) & (a_{[9, 15]}, b_{[9, 15]}) &= (0.5,0.0)\\
(a_{[6, 7]}, b_{[6, 7]}) &= (0.5,0.0) & (a_{[3, 6]}, b_{[3, 6]}) &= (0.0,0.0) & (a_{[10, 12]}, b_{[10, 12]}) &= (0.0,0.5)\\
(a_{[7, 8]}, b_{[7, 8]}) &= (0.5,0.5) & (a_{[3, 13]}, b_{[3, 13]}) &= (0.0,0.5) & (a_{[10, 15]}, b_{[10, 15]}) &= (0.5,0.0)\\
(a_{[8, 9]}, b_{[8, 9]}) &= (0.5,0.5) & (a_{[4, 6]}, b_{[4, 6]}) &= (0.0,0.0) & (a_{[13, 14]}, b_{[13, 14]}) &= (0.0,1.0)\\
& & & & (a_{[15, 16]}, b_{[15, 16]}) &= (1.0,0.0)
\end{align*}

%
%

\eproof

\section{Binomial Edge Ideals of K\"onig type}

\label{sec_ideals_of_Konig_type}

It is a well-known fact that bipartite graphs satisfy the K\"{o}nig--Egervary property, i.e. the size of a minimal vertex cover equals the size of a maximal matching. For (monomial) edge ideals the K\"onig-Egervary property translates into the statement that the monomial grade of the ideal is equal to its height. In \cite{herzog2022graded} the authors introduce the notion of graded ideals of K\"{o}nig type as a way to extend the K\"onig--Egervary property to other homogeneous ideals of a polynomial ring. Because in this section we are only concerned with binomial edge ideals of K\"onig type, we (by abuse of notation) define a graph $G$ to be of K\"onig type if it satisfies the following characterization:

\begin{defn}[\cite{herzog2022graded}  Theorem 3.5]
\label{defn_Konig_type}
A graph $G$ is said to be of {\it K\"onig type} if $r(G) = \hgt(J_{G})$.
\end{defn}

In this section we prove that every graph of K\"onig type satisfies Conjecture \ref{conjectue_fpt_equal_LP}, we give a new combinatorial characterization of the K\"onig type property, and we use this characterization to prove that every tree is of K\"onig type, and hence give a second proof that trees satisfy Conjecture \ref{conjectue_fpt_equal_LP}. Lastly, we give an example of a bipartite graph which is not of K\"onig type.

\bprop
\label{prop_formula_ft_Konig_type}
If $G$ is of K\"onig type, then $\ft(J_{G}) = r(G)$. In particular, Conjecture \ref{conjectue_fpt_equal_LP} holds for every graph $G$ that is of K\"onig type.
\eprop

\bproof
Follows immediately from Equation \eqref{eqn_chain_of_relations} and Definition \ref{defn_Konig_type}.
\eproof

The following lemma gives a combinatorial characterization for a graph to be of K\"onig type. We use this characterization to show that trees are of K\"onig type.

\blem
\label{lem_characterization_Konig_type}
Let $G$ be a graph. Then, $G$ is of K\"onig type if and only if there exists a semi-path $P$ and a choice of distinct vertices $T \subset V(G)$, possibly empty, with the property that
\begin{enumerate}
\item  $r(G) = \card{E(P)}$, \label{item_konig_1}
\item for all $v \in T$, $v$ has degree two in $P$, \label{item_konig_2}
\item for all $v,w \in T$, $\{v,w\}$ is not an edge of $P$, \label{item_konig_3}
\item the number of connected components and the number of vertices of $\widehat{G \setminus T}$ equals the number of connected components and the number of vertices of $\widehat{P \setminus T}$, respectively. \label{item_konig_4}
\end{enumerate}
\elem

\bproof
We prove the backward implication. Suppose that a semi-path $P$ of $G$ and a choice of distinct vertices $T \subset V(P)$ have been chosen satisfying properties \eqref{item_konig_1} through \eqref{item_konig_4}. It is enough to show that $\hgt(J_{G}) \leq r(G)$.

Let $c := c_{G}(T)$ be the number of connected components of $\widehat{G \setminus T}$ and $G_{1},\ldots,G_{c}$ its connected components. Let $Q := (\{X_{i},Y_{i}\}_{i\in T}, J_{\tilde{G}_{1}},\ldots,J_{\tilde{G}_{c}})$ the ideal as defined in Proposition \ref{prop_height_min_primes_bin_edge_ideal}. Then, 
\begin{align*}
\hgt(Q) &= 2\cdot \card{T} + \sum_{i=1}^{c} (\card{V(G_{i})}-1)
\end{align*}
Let $P_{i}$ denote the subgraph of $G_{i}$ given by $G_{i} \cap P$ for $1 \leq i \leq c$. Property \eqref{item_konig_4} implies that $P_{i}$ is a Hamiltonian path of $G_{i}$, and in particular that 
\begin{align*}
\card{E(P_{i})} = \card{V(G_{i})} -1.
\end{align*}
Let $E^{'}$ denote the set of edges of $P$ which are incident to a vertex of $T$. By assumptions \eqref{item_konig_2} and \eqref{item_konig_3}, $\card{E^{'}} = 2\cdot \card{T}$. The edges of $P$ can be partitioned into those which belong to some $P_{i}$ and those which belong to $E^{'}$. We compute that
\begin{align*}
r(G) 
&= \card{E(P)} \phantom{aaaaaaaaaaaaaaaaaaaa} (\text{by \eqref{item_konig_1}}) \\
&= \card{E^{'}} + \sum_{i=1}^{c} \card{E(P_{i})} \\
&= 2\cdot \card{T} + \sum_{i=1}^{c} (\card{V(G_{i})}-1) \\
&= \hgt(Q).
\end{align*}
This proves that $\hgt J_{G} \leq \hgt Q = r(G)$.

Next, suppose that $G$ is of K\"onig type. Let $P$ be a semi-path realizing $r(G)$. Consider the vector $\mbf{a}  = (a_{e})_{e \in E(G)}$ given by 
\begin{align*}
a_{e} = 
\begin{cases}
1, & e \in E(P) \\
0, & \text{otherwise}
\end{cases}
\end{align*}
Then, $\mbf{a} := (a_{e})$ is an optimal solution of the linear program \eqref{alg_path_packing_int_program}. Let $(\mbf{b},\mbf{c}) := (\{b_{v}\}_{v\in V}, \{c_{U}\}_{U \in \A})$ be an optimal solution with values in $\Z$ of the linear program \eqref{alg_dual_path_packing_int_program} having the form guaranteed by Lemma \ref{lem_nice_form_optimal_soln_dual_path_packing}, and let $T$ and $U_{1},\ldots, U_{c}$ be as in Lemma \ref{lem_nice_form_optimal_soln_dual_path_packing}.

Let $A$ be the below $\left(2^{\card{V}}-1\right) \times \card{E}$ matrix, $\mathbf{X}$ the vector of indeterminates of length $\card{E}$, $\mbf{1}$ vector of $1$'s of length $\card{E}$, and $\mbf{d}$ the vector of length $2^{\card{V}-1}$ realizing the linear program \eqref{alg_path_packing_int_program}, i.e. algorithm \eqref{alg_path_packing_int_program} can be expressed as
\begin{equation}
\begin{split}
\maximize \hspace{.3em} & Z_{\PP,G} := \mbf{1}^{\tr} \, \mathbf{X} \\
\subject & \phantom{ } \\
 & \mathbf{X} \in \k_{\geq 0} \\
 & A \cdot \mbf{X} \leq \mbf{d}
\end{split}
\end{equation}
Let $A_{i}$ denote the $i$-th column of $A$, and let $a_{j}^{\tr}$ denote the $j$-th row of $A$.

Because $G$ is of K\"onig type the cost of these optimal solutions are equal, and complementary slackness, Theorem \ref{thm_comp_slackness}, holds between these pair of solutions.

We prove part \eqref{item_konig_2}. Choose $v \in T$, equivalently, $b_{v} = 1$. Let $a_{v}^{\tr}$ denote the row of $A$ corresponding to vertex $v$. By complementary slackness Theorem \ref{thm_comp_slackness}.\eqref{item_comp_slackness_2} we have that
\begin{align*}
2 = a_{v}^{\tr} \mbf{a} = \sum_{e \text{ incident to } v} a_{e}.
\end{align*}
This shows that there are two edges of $P$ incident with $v$.

We prove part \eqref{item_konig_3}. Let $v \in T$, i.e. $b_{v} = 1$, and $e = \{v,w\}$ an edge of $P$, then we want to show that $w \notin T$. Let $A_{e}$ denote column of $A$ corresponding to the edge $e$. Theorem \ref{thm_comp_slackness}.\eqref{item_comp_slackness_1} implies that
\begin{equation}
\label{eqn_pf_of_lem_part_three}
1 = (\mbf{b},\mbf{c})^{\tr} A_{e} = b_{v} + b_{w} + \sum_{\substack{U \in \A \\e\in G[U]}} c_{U}
\end{equation}
Since $b_{v} = 1$, Equation \eqref{eqn_pf_of_lem_part_three} implies that $b_{w} = 0$. Hence $w \notin T$ by Lemma \ref{lem_nice_form_optimal_soln_dual_path_packing}.

We prove part \eqref{item_konig_4}. Fix $1 \leq i \leq c$. Let $a_{U_{i}}^{\tr}$ denote the row of $A$ corresponding to $U_{i}$. By the construction of $\mbf{c}$, $c_{U_{i}} = 1$. Theorem \ref{thm_comp_slackness}.\eqref{item_comp_slackness_2} implies that 
\begin{align*}
\card{U_{i}}-1 = a_{U_{i}}^{\tr} \cdot \mbf{a} = \card{E(P) \cap E(U_{i})}.
\end{align*}
This shows that $P$ restricted to $U_{i}$ induces a Hamiltonian path for every $1 \leq i \leq c$ and \eqref{item_konig_4} follows.
\eproof

\subsection{Trees are of K\"onig Type}

We use Lemma \ref{lem_characterization_Konig_type} to show that trees are of K\"onig type. The proof idea is to pick a maximal semi-path and to recursively delete vertices incident with edges not belonging to the semi-path so at the end we are left with no edges between connected components of the semi-path and no two deleted vertices are adjacent.

We set up notation below. Throughout this subsection let $G$ be a tree. Let $P_{1},\ldots,P_{\ell}$ be disjoint paths in $G$ with $r(G) = \sum_{i=1}^{\ell} \card{E(P_{i})}$. Put 
\begin{align*}
E_{L} := \bigg\{ \vphantom{y^{2}} \{v,w\} \in E(G) \mid 
& v \in P_{i} \text{ and } w \in P_{j} \text{ for some } 1 \leq i < j \leq \ell \\
& \text{ or } \, v\in P_{i} \text{ for some } 1\leq i \leq \ell \text{ and } w\notin P_{j} \text{ for any } 1 \leq j \leq \ell \bigg\}.
\end{align*}
Put $V_{L} := \{V(e) \mid e \in E_{L}\}$ and
\begin{align*}
B_{v} &:= \{ v \in V_{L} \mid v \text{ is an initial or terminal vertex of some } P_{i} \text{ or } v\notin P_{i} \text{ for any } 1\leq i \leq \ell\}.
\end{align*}

\bdefn
\label{defn_subtree_based_at_v}
Let $v \in V(G)$. A {\it subtree of $G$ based at $v$}, denoted $T_{v}$, is the subtree of $G$ constructed as follows:
\begin{align*}
T_{v,k} &:= \varnothing, k \leq 0 \\
T_{v,1} &:= \{v\} \\
T_{v,2k} &:= \{v \in V_{L} \mid \exists w \in T_{2k-1} \text{ with } \{v,w\} \in E_{L}\}, k \geq 1 \\
T_{v,2k+1} &:= \{v \in V(P_{i}) \mid \exists w \in T_{2k} \text{ with } \{v,w\} \in E(P_{i}) \text{ for some } 1 \leq i \leq \ell\}, k \geq 1.
\end{align*}
Then, $T_{v}$ is the induced subgraph of $G$ on the set of vertices $\bigcup_{k\in \Z} T_{k}$.
\edefn

\blem
Let $T_{v}$ be the subtree of $G$ based at $v$. For $N \neq N^{'}$, $T_{v,N} \cap T_{v,N^{'}} = \varnothing$.
\elem

\bproof
Clear because $T_{v,N} \subset \{w \in V(G) \mid d(v,w) = N\}$, and the later sets are disjoint for $N \neq N^{'}$ since $G$ is a tree.
\eproof

The following definition and lemma can be viewed as analogues of the the statement that there are no augmenting paths with respect to a choice of maximal matching in a bipartite graph.

\bdefn
\label{defn_alt_path}
A $(v,w)$-alternating path \textit{with respect to $P$} (respectively a $v$-alternating path \textit{with respect to $P$}) is a set of distinct vertices $v_{1},v_{2},\ldots,v_{2k-1},v_{2k}$ where $1 \leq k \in \Z$ satisfying
\begin{enumerate}
\item $\{v_{2i-1},v_{2i}\} \in E_{L}$ for $1 \leq i \leq k$, \label{item_alt_path_1}
\item $\{v_{2i},v_{2i+1}\} \in E(P_{j_{i}})$ for $1 \leq i \leq k-1$ and for some $1 \leq j_{i} \leq \ell$, \label{item_alt_path_2}
\item $v_{1} = v$ and $v_{2k} = w$ (respectively, $v_{1} = v$). \label{item_alt_path_3}
\end{enumerate}
\edefn

\blem
\label{lem_augmenting_path_restriction}
Let $G$ be a graph. Then, $G$ does not have an alternating path $v_{1}$, $v_{2}, \ldots,$ $v_{2k-1}$, $v_{2k}$ with both $v_{1} \in B_{v}$ and $v_{2k} \in B_{v}$.
\elem

\bproof
If such an alternating path were to exist, then we could construct a semi-path having more edges than $P$ by adding to $P$ the edges belonging to \eqref{item_alt_path_1} of Definition \ref{defn_alt_path} and deleting from $P$ the edges belonging to \eqref{item_alt_path_2} of Definition \ref{defn_alt_path}.
\eproof

\balg
\label{alg_showing_trees_Konig_type}
$\newline$
\textbf{Input}: Semi-path $P$ of a tree $G$ $\newline$
\textbf{Output}: $T \subseteq \bigcup_{i=1}^{\ell} V(P_{i})$ satisfying conditions \eqref{item_konig_2} - \eqref{item_konig_3} in Lemma \ref{lem_characterization_Konig_type}, and every edge of $E_{L}$ contains a vertex belonging to $T$.

\begin{enumerate}
\item Initialize $i = 0$ and $T_{0,k} = \varnothing$ for $k \in \Z$.
\item While $B_{v} \not\subseteq \bigcup\limits_{n=0}^{i} \bigcup\limits_{k\geq 1} T_{n,k}$, repeat the following steps: \label{item_alg_tree_is_konig_part_2}
\begin{enumerate}
\item Increment $i$,
\item Pick a vertex $v_{i} \in B_{v} \setminus \bigcup\limits_{n=0}^{i-1} \bigcup\limits_{k\geq 1} T_{n,k}$.
\item Let $T_{i}$ denote the subtree of $G$ based at $v_{i}$, and for $k \in \Z$, let $T_{i,k}$ denote the subsets of $T_{i}$ as constructed in Definition \ref{defn_subtree_based_at_v}.
\end{enumerate}
\item Put $S := V_{L} \setminus \bigcup\limits_{n=0}^{i} \bigcup\limits_{k\geq 1} T_{n,k}$. \label{item_alg_tree_is_konig_part_3}
\item While $S \neq \varnothing$, repeat the following steps: 
\begin{enumerate}
\item Increment $i$,
\item Choose $1 \leq j \leq \ell$ minimally so that there exists a vertex $v_{i} \in V(P_{j}) \cap S$,
\item Let $T_{i}$ denote the subtree of $G$ based at $v_{i}$, and for $k \in \Z$, let $T_{i,k}$ denote the subsets of $T_{i}$ as constructed in Definition \ref{defn_subtree_based_at_v},
\item Go back to part \eqref{item_alg_tree_is_konig_part_3} of Algorithm \ref{alg_showing_trees_Konig_type}.
\end{enumerate}
\item Put $N := i$. 
\item For $1 \leq i \leq N$, define $T^{\even}_{i} := \bigcup_{k\in \Z} T_{i,2k}$ and $T^{\odd}_{i} := \bigcup_{k\in \Z} T_{i,2k+1}$.
\item Put $T := \bigcup_{i=1}^{N} T_{i}^{\even}$.
\end{enumerate}
\ealg

\bthm
\label{thm_alg_constructing_T_correct}
Algorithm \ref{alg_showing_trees_Konig_type} is correct, and property \eqref{item_konig_4} of Lemma \ref{lem_characterization_Konig_type} is satisfied by the choice of $T$. Consequently, $G$ is of K\"onig type.
\ethm

\bproof
First, we remark that each while-loop of Algorithm 
\ref{alg_showing_trees_Konig_type} terminates because $B_{v}$ and $V_{L}$ are finite, and $\{ \bigcup_{n=0}^{i} \bigcup_{k\geq 1} T_{n,k} \}_{i=1}^{N}$ form a strictly ascending sequence of sets each containing more elements of $B_{v} \cup V_{L}$ than their predecessors.

\begin{claim}
\label{claim_disjnt_sets}
For $i \neq j$, $T^{\even}_{i} \cap T^{\odd}_{j} = \varnothing$.
\end{claim}	

\begin{pfclaim}
Suppose by contradiction that $T^{\even}_{i} \cap T^{\odd}_{j} \neq \varnothing$. Choose $a$ even and $b$ odd satisfying the following
\begin{enumerate}[(i)]
\item there exists $w \in T_{i,a} \cap T_{j,b}$,
\item if $q_{1}$ (respectively $q_{2}$) is the $(v_{i},w)$-alternating path (respectively $(v_{j},w)$-alternating path), then the vertices of $q_{1}$ and $q_{2}$ intersect only at $w$.
\end{enumerate}
Let $w^{'}$ be the vertex of $q_{2}$ adjacent to $w$. Then, there is a $(v_{i},v_{j})$-augmenting path by concatenating $q_{1}$ and $q_{2}$ along the edge $\{w,w^{'}\}$. Without loss of generality, suppose that $i < j$, then $v_{j} \in V(T_{i})$ by construction of $T_{i}$, but this contradicts the choice of $v_{j} \in V_{L} \setminus \bigcup_{n=1}^{j-1} \bigcup_{k\geq 1} T_{n,k}$ which proves the claim.
\end{pfclaim}

\ul{$T \subset \bigcup_{i=1}^{\ell} V(P_{i})$}: Suppose by contradiction that there exists $w \in T \setminus \bigcup_{i=1}^{\ell} V(P_{i})$. Since $w \in B_{v}$ and $w \in T_{i}^{\even}$ for some $1 \leq i \leq N$, it must be the case that $v_{i} \in B_{v}$ by part \ref{item_alg_tree_is_konig_part_2} of Algorithm \ref{alg_showing_trees_Konig_type}. This contradicts Lemma \ref{lem_augmenting_path_restriction}.

\ul{$T$ satisfies \eqref{item_konig_2} of Lemma \ref{lem_characterization_Konig_type}}: $T$ does not contain any initial nor terminal points of any $P_{i}$ by a similar argument to the previous paragraph.

\ul{$T$ satisfies \eqref{item_konig_3} of Lemma \ref{lem_characterization_Konig_type}}: If not, we would contradict $T_{i}^{\even} \cap T_{j}^{\odd} = \varnothing$.

\ul{$T$ satisfies \eqref{item_konig_4} of Lemma \ref{lem_characterization_Konig_type}}: The construction of $T$ shows that for each edge $\{v,w\} \in E_{L}$ that
\begin{enumerate}[(i)]
\item $w \in T$ if $v \notin \bigcup_{i=1}^{\ell} P_{i}$, and that
\item $v \in T$ or $w \in T$ if $\{v,w\} \in E_{L}$.
\end{enumerate}
Thus, as sets $\wh{G\setminus T}$ and $\wh{P \setminus T}$ are equal, and hence condition \eqref{item_konig_4} is trivially satisfied. Hence by Lemma \ref{lem_characterization_Konig_type} $G$ is of K\"onig type.
\eproof

\subsection{Example of Bipartite Graph Not of K\"onig Type}

The following example shows that there exist bipartite graphs which are not of K\"onig type.

\def\putAbove{-.35}
\def\edgeDistBip{.75}

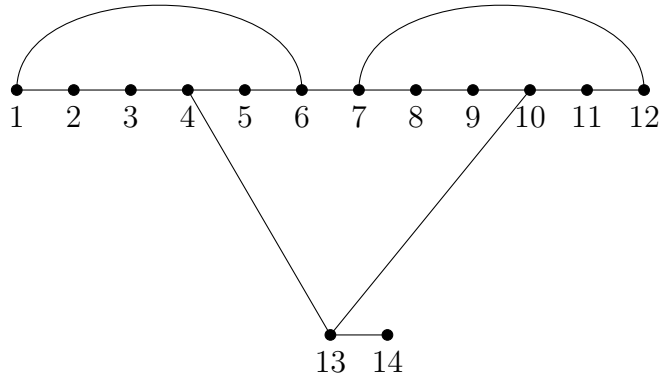
\begin{figure}[h]
\begin{center}
\begin{tikzpicture} 
\filldraw[black] (-5.5*\edgeDistBip,0) circle (2pt) node at (-5.5*\edgeDistBip,\putAbove) {$1$};
\filldraw[black] (-4.5*\edgeDistBip,0) circle (2pt) node at (-4.5*\edgeDistBip,\putAbove) {$2$};
\filldraw[black] (-3.5*\edgeDistBip,0) circle (2pt) node at (-3.5*\edgeDistBip,\putAbove) {$3$};
\filldraw[black] (-2.5*\edgeDistBip,0) circle (2pt) node at (-2.5*\edgeDistBip,\putAbove) {$4$};
\filldraw[black] (-1.5*\edgeDistBip,0) circle (2pt) node at (-1.5*\edgeDistBip,\putAbove) {$5$};
\filldraw[black] (-.5*\edgeDistBip,0) circle (2pt) node at (-.5*\edgeDistBip,\putAbove) {$6$};

\filldraw[black] (.5*\edgeDistBip,0) circle (2pt) node at (.5*\edgeDistBip,\putAbove) {$7$};
\filldraw[black] (1.5*\edgeDistBip,0) circle (2pt) node at (1.5*\edgeDistBip,\putAbove) {$8$};
\filldraw[black] (2.5*\edgeDistBip,0) circle (2pt) node at (2.5*\edgeDistBip,\putAbove) {$9$};
\filldraw[black] (3.5*\edgeDistBip,0) circle (2pt) node at (3.5*\edgeDistBip,\putAbove) {$10$};
\filldraw[black] (4.5*\edgeDistBip,0) circle (2pt) node at (4.5*\edgeDistBip,\putAbove) {$11$};
\filldraw[black] (5.5*\edgeDistBip,0) circle (2pt) node at (5.5*\edgeDistBip,\putAbove) {$12$};

\filldraw[black] (0*\edgeDistBip,-4.33012701892*\edgeDistBip) circle (2pt) node at (0*\edgeDistBip,-4.33012701892*\edgeDistBip+\putAbove) {$13$};
\filldraw[black] (1*\edgeDistBip,-4.33012701892*\edgeDistBip) circle (2pt) node at (1*\edgeDistBip,-4.33012701892*\edgeDistBip+\putAbove) {$14$};

\draw (-5.5*\edgeDistBip,0) to (-4.5*\edgeDistBip,0);
\draw (-4.5*\edgeDistBip,0) to (-3.5*\edgeDistBip,0);
\draw (-3.5*\edgeDistBip,0) to (-2.5*\edgeDistBip,0);
\draw (-2.5*\edgeDistBip,0) to (-1.5*\edgeDistBip,0);
\draw (-1.5*\edgeDistBip,0) to (1.5*\edgeDistBip,0);
\draw (1.5*\edgeDistBip,0) to (2.5*\edgeDistBip,0);
\draw (2.5*\edgeDistBip,0) to (3.5*\edgeDistBip,0);
\draw (3.5*\edgeDistBip,0) to (4.5*\edgeDistBip,0);
\draw (4.5*\edgeDistBip,0) to (5.5*\edgeDistBip,0);

\draw (-2.5*\edgeDistBip,0) to (0,-4.33012701892*\edgeDistBip);
\draw (3.5*\edgeDistBip,0) to (0,-4.33012701892*\edgeDistBip);
\draw (0,-4.33012701892*\edgeDistBip) to (1*\edgeDistBip,-4.33012701892*\edgeDistBip);

\draw (-5.5*\edgeDistBip,0) .. controls (-5.5*\edgeDistBip,2*\edgeDistBip) and (-.5*\edgeDistBip,2*\edgeDistBip) .. (-.5*\edgeDistBip,0);

\draw (.5*\edgeDistBip,0) .. controls (.5*\edgeDistBip,2*\edgeDistBip) and (5.5*\edgeDistBip,2*\edgeDistBip) .. (5.5*\edgeDistBip,0);
\end{tikzpicture}
\end{center}
\caption{Non-K\"onig type Non-Unmixed Bipartite Graph}
\label{fig_non_konig_bipartite_graph}
\end{figure}

\bex
\label{ex_non_konig_bipartite_graph}
For this example let $G$ be the graph in Figure \ref{fig_non_konig_bipartite_graph}.

One can check that $r(G) = 12$, $\LP_{G,\Q} = 12.5$, and $\hgt J_{G} = 13$. One realization of $\LP_{G,\Q}$ is via the following assignment
%
%
%
%
%
\begin{align*}
(a_{[1, 2]}, b_{[1, 2]}) &= (0.5,0.5)  & (a_{[6, 7]}, b_{[6, 7]}) &= (0.0,0.5)  & (a_{[11, 12]}, b_{[11, 12]}) &= (1.0,0.0) \\
(a_{[2, 3]}, b_{[2, 3]}) &= (0.5,0.5) & (a_{[7, 8]}, b_{[7, 8]}) &= (0.5,0.0) & (a_{[13, 14]}, b_{[13, 14]}) &= (0.5,0.5) \\
(a_{[3, 4]}, b_{[3, 4]}) &= (0.5,0.5) & (a_{[8, 9]}, b_{[8, 9]}) &= (0.5,0.5) & (a_{[1, 6]}, b_{[1, 6]}) &= (0.0,0.5) \\
(a_{[4, 5]}, b_{[4, 5]}) &= (0.0,0.5) & (a_{[9, 10]}, b_{[9, 10]}) &= (0.5,0.5) & (a_{[7, 12]}, b_{[7, 12]}) &= (0.0,1.0) \\
(a_{[5, 6]}, b_{[5, 6]}) &= (0.5,0.5) & (a_{[10, 11]}, b_{[10, 11]}) &= (0.5,0.0) & (a_{[4, 13]}, b_{[4, 13]}) &= (0.5,0.0) \\
& & & &  (a_{[10, 13]}, b_{[10, 13]}) &= (0.0,0.5)
\end{align*}
This graph is not unmixed. For example, $\{1,3\}$ is a cut set of $G$ and the corresponding prime ideal has height $14$. This raises
\eex

\bques
\label{ques_unmixed_bipartite_implies_Konig_type}
Is every unmixed bipartite graph of K\"onig type?
\eques

Bipartite graphs having unmixed binomial edge ideal were studied in \cite{bolognini2018binomial}. Bipartite graphs with Cohen--Macaulay binomial edge ideal were characterized in {\cite[Theorem 6.1]{bolognini2018binomial}}. As a consequence of their characterization, Cohen--Macaulay bipartite graphs are traceable. Herzog et al. proved that traceable implies K\"onig type {\cite[Proposition 3.6.(a)]{herzog2022graded}}, and thus in particular Cohen--Macaulay bipartite graphs are of K\"onig type.

\section{$\LP_{G,\Q}$ and Traceable Graphs}

\label{sec_ft_and_Ham_paths}

In this section, we consider the relationship between maximality of $\LP_{G,\Q}$ of a graph and the graph being traceable. Recall that a graph is called traceable if it has a Hamiltonian path, Definition \ref{defn_traceable_graphs}.

\bprop
\label{prop_ft_traceable_graph}
If $G$ is a traceable graph on $n$ vertices, then $\ft(J_{G}) = \LP_{G,\Q} = n-1$.
\eprop

\bproof
We have the graph containments $P_{n} \subset G \subset K_{n}$. The result follows from monotonicity of $\LP_{G,\k}$ and $\ft(J_{P_{n}}) = \LP_{K_{n},\Q} = n-1$.
\eproof

\brem
The above proof shows that for any graph $G$ on $n$ vertices, $\ft(J_{G}) \leq n-1$.
\erem

Example \ref{ex_LP=n-1_non-traceable} shows that there exists a graph $G$ with $\LP_{G,\Q} = \card{V} -1$ yet $G$ is non-traceable. Computations in Macaulay2 \cite{M2} suggest that $\ft(J_{G})$ may agree with $\card{V}-1$ independent of characteristic in Example \ref{ex_LP=n-1_non-traceable}.


\bex
\label{ex_LP=n-1_non-traceable}
Throughout this example let $G$ denote the graph in Figure \ref{fig_non_traceable_graph}. Observe that $G$ is not traceable. 

\def\putBelow{-.5}
\def\edgeDist{1.75}

\begin{figure}[h]

\begin{center}
\begin{tikzpicture} 
\filldraw[black] (-3*\edgeDist,0) circle (2pt) node at (-3*\edgeDist,\putBelow) {1};
\filldraw[black] (-2*\edgeDist,0) circle (2pt) node at (-2*\edgeDist,\putBelow) {2};
\filldraw[black] (-1*\edgeDist,0) circle (2pt) node at (-1*\edgeDist,\putBelow) {4};
\filldraw[black] (0,0) circle (2pt) node at (0,\putBelow) {6};
\filldraw[black] (1*\edgeDist,0) circle (2pt) node at (1*\edgeDist,\putBelow) {8};
\filldraw[black] (2*\edgeDist,0) circle (2pt) node at (2*\edgeDist,\putBelow) {9};
\filldraw[black] (3*\edgeDist,0) circle (2pt) node at (3*\edgeDist,\putBelow) {10};
\filldraw[black] (-3/2*\edgeDist, 
0.866025403784*\edgeDist) circle (2pt) node at (-3/2*\edgeDist, 
0.866025403784*\edgeDist+.5) {3};
\filldraw[black] (3/2*\edgeDist, 
0.866025403784*\edgeDist) circle (2pt) node at (3/2*\edgeDist, 
0.866025403784*\edgeDist+.5) {7};
\filldraw[black] (0,1.75*\edgeDist) circle (2pt) node at (0,1.75*\edgeDist+0.5) {5};

\draw (-3*\edgeDist,0) to (-2*\edgeDist,0);
\draw (-2*\edgeDist,0) to (-1*\edgeDist,0);
\draw (-1*\edgeDist,0) to (0,0);
\draw (-0,0) to (1*\edgeDist,0);
\draw (1*\edgeDist,0) to (2*\edgeDist,0);
\draw (2*\edgeDist,0) to (3*\edgeDist,0);
\draw (-2*\edgeDist,0) to (-3/2*\edgeDist, 
0.866025403784*\edgeDist);
\draw (-1*\edgeDist,0) to (-3/2*\edgeDist, 
0.866025403784*\edgeDist);
\draw (2*\edgeDist,0) to (3/2*\edgeDist, 
0.866025403784*\edgeDist);
\draw (1*\edgeDist,0) to (3/2*\edgeDist, 
0.866025403784*\edgeDist);
\draw (0,1.75*\edgeDist) to (-3/2*\edgeDist, 
0.866025403784*\edgeDist);
\draw (0,1.75*\edgeDist) to (3/2*\edgeDist, 
0.866025403784*\edgeDist);

\end{tikzpicture}
\end{center}
\caption{Non-traceable graph}
\label{fig_non_traceable_graph}
\end{figure}
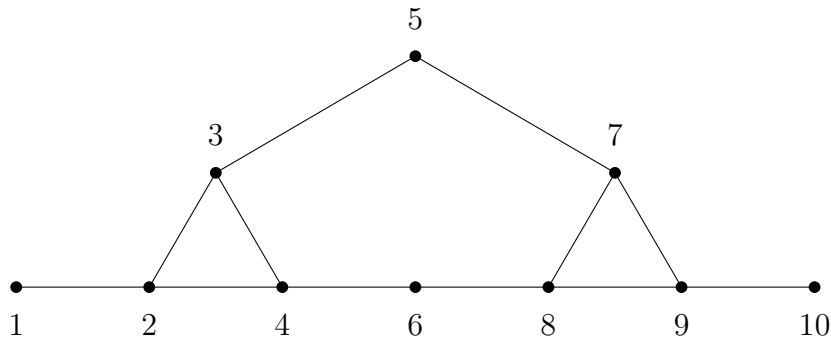

Observe that deleting edges $\{7,8\}$ and $\{8,9\}$ produces a triangle path as a subgraph and that this triangle path has F-threshold equal to $8.5$ by Proposition \ref{prop_formula_for_fpt_of_triangle_path}. It follows from Proposition \ref{lem_ft_monotonic} that $8.5 \leq \ft(J_{G})$. However, it can be computed that $\LP_{G,\Q} = 9$ and one such realization is as follows
\begin{align*}
(a_{[1,2]},b_{[1,2]}) &= (1.0,0.0) & (a_{[3,5]},b_{[3,5]}) &= (1.0,0.0) & (a_{[7,8]},b_{[7,8]}) &= (0.5,0.0) \\
(a_{[2,3]},b_{[2,3]}) &= (0.5,0.0) & (a_{[4,6]},b_{[4,6]}) &= (0.5,0.5) & (a_{[7,9]},b_{[7,9]}) &= (0.5,0.0) \\
(a_{[2,4]},b_{[2,4]}) &= (0.5,0.0) & (a_{[5,7]},b_{[5,7]}) &= (1.0,0.0) & (a_{[8,9]},b_{[8,9]}) &= (0.5,0.0) \\
(a_{[3,4]},b_{[3,4]}) &= (0.0,0.5) & (a_{[6,8]},b_{[6,8]}) &= (0.5,0.5) & (a_{[9,10]},b_{[9,10]}) &= (1.0,0.0)
\end{align*}
%
%
%
%
%
%
This bilabeling corresponds to the monomial 
\begin{align*}
\prod_{i=0}^{8} X_{i} \cdot \prod_{i=1}^{9} Y_{i}
\end{align*}
belonging to the polynomial
\begin{align*}
\big( f_{01} f_{13}^{1/2} f_{12}^{1/2} f_{23}^{1/2} f_{24} f_{35} f_{46} f_{57} f_{67}^{1/2} f_{68}^{1/2} f_{78}^{1/2} f_{89} \big)^{p^{e}-1} \in J_{G}^{9(p^{e}-1)}
\end{align*}
which has coefficient given by the expression
\begin{align*}
\sum_{a=0}^{N/2} \bigg( \binom{N/2}{a} \binom{N}{N/2+a} \binom{N}{a} \bigg)^{2}
\end{align*}
where $N := p^{e}-1$. For small values of $p$ and $e$ it can be checked that this coefficient does not vanish mod $p$. This suggests that $\ft(J_{G}) = 9$ for any prime $p$.
\eex

The following proposition shows that under additional hypotheses on $G$ that a converse of Proposition \ref{prop_ft_traceable_graph} holds.

\bprop
\label{three_verts_val_one}
If $G$ is a connected graph on $n+3$ vertices with at least three vertices having degree one (in particular $G$ is not traceable), then $\LP_{G,\Q} \leq n + 3/2 < n+2$.
\eprop

\bproof
Observe that $G$ is a subgraph of a complete graph $H$ on $n$ vertices with three whiskers. Suppose that the whiskers of $H$ do not share a common vertex, then by Corollary  \ref{cor_ft_two_block_graph_equals_ft_triangle_path_packing}, $\LP_{H,\Q} = n + \frac{3}{2}$. If two of the three whiskers or all three of the whiskers of $H$ share a common vertex, then we can list the two block approximations of $H$ and show that $\LP_{H,\Q} \leq n+1$ by Theorem \ref{thm_formula_ft_block_graph}. Since $\LP_{G,\Q} \leq \LP_{H,\Q}$ the result follows.
\eproof

It can be shown that the graph in Example \ref{ex_LP=n-1_non-traceable} is not unmixed. This raises

\bques
\label{ques_unmixed_plus_LP=n-1_implies_trac}
If $G$ is an unmixed, connected graph and $\LP_{G,\Q} = \card{V(G)}-1$, then is $G$ traceable?
\eques

Note that $G$ being unmixed and connected implies that $\LP_{G,\Z}^{\ast} = n-1$. If Question \ref{ques_unmixed_plus_LP=n-1_implies_trac} were true, this would strengthen a result of \cite[3.6.(b)]{herzog2022graded} that an unmixed graph of K\"onig type is traceable. In terms of the notation of this paper, their statement can be rephrased as: if $G$ is unmixed, connected, and $\LP_{G,\Z} = \LP_{G,\Z}^{\ast}$, then $G$ is traceable. Question \ref{ques_unmixed_plus_LP=n-1_implies_trac} can be rephrased as: if $G$ is  unmixed, connected, and $\LP_{G,\Q} = \LP_{G,\Z}^{\ast}$, then is $G$ traceable?

 \bibliographystyle{amsalpha}
  \bibliography{Bibliography_Invariants_of_Binomial_Edge_Ideals_via_Linear_Programs}

\end{document}